\documentclass[review]{elsarticle}

\usepackage{amssymb}
\usepackage{amsmath}
\usepackage{amscd}
\usepackage{amsthm}
\usepackage{graphicx}
\usepackage{mathrsfs}
\usepackage{amscd}

\biboptions{sort&compress}

\newenvironment{pf}[1][Proof.]{\noindent\textbf{#1 }\\ }{\hfill$\square$\\}

\newtheorem{thm}{Theorem}[section]
\newtheorem{lem}[thm]{Lemma}
\newtheorem{prop}[thm]{Proposition}
\newtheorem{cor}[thm]{Corollary}
\newtheorem{defn}[thm]{Definition}

\numberwithin{equation}{section}

\newcommand{\sA}{\mathscr{A}}
\newcommand{\sB}{\mathscr{B}}
\newcommand{\sC}{\mathscr{C}}
\newcommand{\sD}{\mathscr{D}}
\newcommand{\sF}{\mathscr{F}}
\newcommand{\sY}{\mathscr{Y}}

\newcommand{\sX}{\mathcal{X}}

\newcommand{\Pf}{\mathcal{P}_{\hspace{-0.2mm}\mathrm{fin}}   }

\newcommand{\hdelta}{{\widehat{\delta}}}
\newcommand{\deltap}[1]{{\delta\left(#1\right)}}

\newcommand{\diam}{\mathrm{diam}}

\newcommand{\conv}[1] {\mathrm{conv}(#1)}

\journal{Advances in Mathematics}

\begin{document}

\begin{frontmatter}

\title{Hyperconvexity and Tight-Span Theory for Diversities}
\author[otago]{David Bryant\corref{cor1}}
\ead{david.bryant@otago.ac.nz}

\author[sf]{Paul F. Tupper}
\ead{pft3@math.sfu.ca}

\cortext[cor1]{Corresponding author}

\address[otago]{Dept. of Mathematics and Statistics, University of Otago. PO Box 56 Dunedin 9054, New Zealand. Ph (64)34797889. Fax (64)34798427 }
\address[sf]{Dept. of Mathematics, Simon Fraser University. 8888 University Drive, Burnaby, British Columbia V5A 1S6, Canada.  Ph (778)7828636. Fax (778)7824947}

\begin{abstract}
The tight span, or injective envelope, is an elegant and useful construction that takes a metric space and returns the smallest hyperconvex space into which it can be embedded. The concept has stimulated a large body of theory and has applications to metric classification and data visualisation. 
Here we introduce a generalisation of metrics, called diversities, and demonstrate that the rich theory associated to metric tight spans and hyperconvexity extends to a seemingly richer theory of diversity tight spans and hyperconvexity.
\end{abstract}

\begin{keyword}
Tight span; Injective hull; Hyperconvex; Diversity; Metric geometry; \end{keyword}

\end{frontmatter}

\section{Introduction}

Hyperconvex metric spaces were defined by Aronszajn and Panitchpakdi in \cite{Aronszajn56} as part of a program to generalise 
 the Hahn-Banach theorem to more general metric spaces (reviewed in \cite{Espinola01}, and below). Isbell \cite{isbell} and Dress \cite{Dress84} showed that, for every metric space, there exists an essentially unique ``minimal'' hyperconvex space into which that space could be embedded, called the {\em tight span} or {\em injective envelope}. Our aim is to show that the notion of hyperconvexity, the tight span, and much of the related theory can be extended beyond metrics to a class of multi-way metrics which we call diversities. 

Recall that a {\em metric space} is a pair $(X,d)$ where $X$ is a set and $d$ is a function from $X \times X$ to   $\Re$ satisfying
\begin{enumerate}
\item[(M1)] $d(a,b) \geq 0$ and $d(a,b) = 0$ if and only if $a=b$.
\item[(M2)] $d(a,c) \leq d(a,b) + d(b,c)$
\end{enumerate}
for all $a,b,c \in X$. 
We define a {\em diversity} to be a pair $(X,\delta)$ where $X$ is a set and $\delta$ is a function from the finite subsets of $X$ to $\Re$ satisfying
\begin{enumerate}
\item[(D1)] $\delta(A) \geq 0$, and $\delta(A) = 0$  if and only if $|A| \leq 1$.
\item[(D2)] If $B \neq \emptyset$ then $\delta(A \cup C) \leq \delta(A \cup B) + \delta(B \cup C)$
\end{enumerate}  
for all finite $A,B,C \subseteq X$.
We prove below that these axioms imply monotonicity:
\begin{enumerate}
\item[(D3)] If $A \subseteq B$ then $\delta(A) \leq \delta(B).$
\end{enumerate}

We will show that tight-span theory adapts elegantly from metric spaces to diversities. The tight span of a metric space $(X,d)$ is formed from the set of point-wise minimal functions $f:X \rightarrow \Re$  such that $f(a_1) + f(a_2) \geq d(a_1,a_2)$ for all $a_1,a_2 \in X$.  Letting $\Pf(X)$ denote the finite subsets of $X$, the tight span $T_X$ of a diversity $(X,\delta)$ is formed from the set of point-wise minimal functions $f:\Pf(X) \rightarrow \Re$  such that 
\[
f(A_1)+f(A_2) + \cdots + f(A_k) \geq \delta(A_1 \cup A_2 \cup \cdots \cup A_k)
\]
 for all finite collections $\{A_1,A_2,\ldots,A_k\} \subseteq \Pf(X)$. The metric tight span is itself a metric space with a canonically defined metric $d_T$; the diversity tight span is itself a diversity $(T_X,\delta_T)$ with a canonically defined function $\delta_T: \Pf(T_X)\rightarrow \Re$. A metric space can be embedded in its tight span; a diversity $(X,\delta)$ can be embedded in its tight span  $(T_X,\delta_T)$. Both constructions have characterisations in terms of injective hulls, and both possess a rich mathematical structure.

The motivation for exploring tight spans of diversities was the success of the metric tight span as a tool for classifying and visualising finite metrics, following the influential paper of Dress \cite{Dress84}. The construction provided the theoretical framework for split decomposition \cite{Bandelt92} and Neighbor-Net \cite{Bryant04}, both implemented in the SplitsTree package  \cite{Huson06} and widely used for visualising phylogenetic data. By looking at diversities, rather than metrics or distances, our hope is to incorporate more information into the analysis and thereby improve inference \cite{Pachter04}.

Dress et al.~\cite{Dress96} coined the term {\em T-theory} for the field of discrete mathematics devoted to the combinatorics of the tight span and related constructions. Sturmfels \cite{Sturmfels05} highlighted T-theory as  one area where problems from biology  have led to substantial new ideas in mathematics. Contributions to T-theory include profound results on optimal graph realisations of metrics \cite{Dress84,Dress06,Huber08};  intriguing connections between the Buneman graph, the tight span and related constructions \cite{Buneman74,Dress96,Dress97,Dress97a,Huber06,Huber08}; links with tropical geometry and hyperdeterminants \cite{Develin04,Huggins08}; classification of finite metrics \cite{Dress84,Sturmfels04}; and properties of the tight span for special classes of metrics \cite{Huber05,Eppstein10}. Hirai \cite{Hirai06} describes an elegant geometric formulation of the tight span. Herrman and Moulton \cite{Herrman11} have recently shown how this geometric framework can be used to study the diversity tight spans which we introduce here.  We believe that there will be diversity analogues for many other metric-space results.

Our use of the term {\em diversity} comes from the appearance of a special case of our definition in work on phylogenetic and ecological diversity \cite{Faith92,Pachter04,Steel05,Minh09}.  However diversities crop up in a broad range of contexts, for example:
\begin{enumerate}
\item 
\noindent {\em Diameter Diversity.}  Let $(X,d)$ be a metric space.  For all $A \in \Pf(X)$ let 
\[
\delta (A) = \max_{a,b \in A} d(a,b) = \diam(A).
\]
Then $(X,\delta)$ is a diversity. 
\item {\em $L_1$ diversity.} For all finite $A \subseteq \Re^n$ define 
\[\delta(A) = \sum_i \max_{a,b} \{|a_i - b_i|:a,b \in A\}.\]
Then $(\Re^n,\delta)$ is a diversity. This result follows directly from the fact that if $(X,\delta_X)$ and $(Y,\delta_Y)$ are diversities and $\delta$ is defined for all finite subsets of $X \times Y$ by 
\[\delta\big(\{(x_1,y_1),\ldots,(x_k,y_k)\}\big) = \delta_X(\{x_1,\ldots,x_k\}) + \delta_Y(\{y_1,\ldots,y_k\})\]
then $(X \times Y,\delta)$ is a diversity.
\item {\em Phylogenetic Diversity.}
Let $T$ be a phylogenetic tree with taxon set $X$. For each finite $A \subseteq X$, let $\delta(A)$ denote the length of the smallest subtree of $T$ connecting taxa in $A$.  Then $(X,\delta)$ is a \emph{phylogenetic diversity}.
\item {\em Length of the Steiner Tree.} 
Let $(X,d)$ be a metric space. For each finite $A \subseteq X$ let $\delta(A)$ denote the minimum length of a Steiner tree within $X$ connecting elements in $A$.   Then $(X,\delta)$ is a diversity.
\item {\em Truncated diversity.} Let $(X,\delta)$ be a diversity. For all $A \in \Pf(X)$ define
\[\delta^{(k)}(A) = \max \{ \delta(B) :|B| \leq k,\,\,B \subseteq A \}.\]  For each $k \geq 2$, $(X,\delta^{(k)})$ is a diversity.
Note that these diversities can be encoded using $O(|X|^k)$ values, an important consideration when designing efficient algorithms.
\end{enumerate}

The generalisation of metrics  to more than two arguments has a long history.    There is an extensive literature on 2-metrics (metrics taking three points as arguments); see \cite{gahler}.  Generalised metrics defined on $n$-tuples for arbitrary $n$ go back at least to  Menger \cite{menger}, who took the volume of an $n$-simplex in Euclidean space as the prototype. Recently various researchers have continued the study of such generalised metrics defined on $n$-tuples; see anal
\cite{deza_rosenberg2,Chepoi_fichet,Warrens}
for examples.  However, as of yet, a satisfactory theory of tight spans has not been developed for these generalisations.  

Dress and Terhalle \cite{Dress98}  developed  tight-span theory for {\em valuated matroids}, which can be viewed as an $n$-dimensional version of a restricted class of metrics. They demonstrated intriguing links with algebraic building theory.  One significant difference is that, for diversities, the tight span consists of functions on $\Pf(X)$ rather than on $X$, as is the case for valuated matroids. 
 
We note that our results  differ from all of  this earlier work because, for a diversity $(X,\delta)$, the function $\delta$ is defined on arbitrary finite subsets of $X$ rather than tuples of a fixed length. In this way, diversities can be compared to valuated $\Delta$-matroids \cite{Dress91}.\\

The structure of this paper is as follows:
In Section~\ref{sec:tightspan} we develop the basic theory of tight spans on diversities. We define the  {\em diversity tight span} $(T_X,\delta_T: \Pf(T_X)\rightarrow \Re)$ of a diversity $(X,\delta)$ and show that every diversity embeds into its diversity tight span. In Section~\ref{sec:hyperconvex} we characterise diversities that are isomorphic to their tight spans. Here, isomorphism is defined in analogy to isometry for metric spaces. These are the {\em hyperconvex} diversities, a direct analogue of hyperconvex metrics. We prove that diversity tight spans, like metric tight spans, are injective, and are formally the injective envelope in the category of diversities.
In Section~\ref{sec:tightlinks} we explore in more detail the direct links between diversity tight spans and metric tight spans. We show that when the diversity equals the diameter diversity (as defined above) the diversity tight span is isomorphic to the diameter diversity of the metric tight span.
In Section~\ref{sec:phylo} we study the tight span of a phylogenetic diversity, and prove that taking the tight span of a phylogenetic diversity recovers the underlying tree in the same way that taking the tight span of an additive metric recovers its underlying tree. This theory is developed for {\em $\Re$-trees}.

Finally, in Section~\ref{sec:Steiner} we examine applications of the theory to the classical Steiner Tree problem. Dress and Kr\"uger \cite{Dress87} defined an {\em abstract} Steiner tree where the internal nodes did not have to sit in the given metric space. They proved that these abstract Steiner trees can be embedded in the tight span. We extend their results to Steiner trees based on diversities, thereby obtaining tight bounds for the classical Steiner tree problem.



\section{The tight span of a diversity} \label{sec:tightspan}

We begin by establishing some basic properties of diversities. Recall that $\Pf(X)$ denotes the set consisting of all finite subsets of the set $X$, and that a diversity is a pair $(X,\delta)$ where $X$ is a set and  the function $\delta \colon \Pf(X) \rightarrow  \Re$ satisfies axioms (D1) and (D2). 

\begin{prop} \label{thm:diversitybasics}
Let $(X,\delta)$ be a diversity.
\begin{enumerate}
\item If $d \colon X \times X \rightarrow \Re$ is defined as $d(x,y)=\delta(\{x,y\})$ then $(X,d)$ is a metric space. We say that $(X,d)$ is the {\em induced metric} of $(X,\delta)$.
\item (D3) holds, that is, for $A,B \in \Pf(X)$, if $A  \subseteq B$ then $\delta(A) \leq \delta(B)$.
\item For $A,B \in \Pf(X)$ if  $A \cap B \neq \emptyset$ then $\delta(A \cup B) \leq \delta(A) + \delta(B)$.
\end{enumerate}
\end{prop}
\begin{pf}
\begin{enumerate}
\item We have  $d(x,y)=0$ if and only if $x=y$ in view of (D1).  Symmetry is clear, and using (D2) we obtain the triangle inequality
\[
d(x,z) = \delta(\{x,z\}) \leq \delta(\{x,y\}) + \delta(\{y,z\}) = d(x,y) + d(y,z).
\]
for all $x,y,z \in X$.
 
\item
First note for any $a \in A$ and $b \in X$ that  by (D2) with $C$ empty
\[
\delta(A) \leq  \delta(A \cup \{b\} )   + \delta(\{b\}) = \delta(A \cup \{b\}). 
\]
The more general result follows by induction.
\item 
Using (D2)
we have
\[ \delta(A \cup B) \leq \delta(A \cup (A \cap B) ) + \delta( B \cup (A \cap B)) = \delta(A) + \delta(B). \]
\end{enumerate}
\end{pf}


We now state the diversity analogue for the metric tight span. 

\begin{defn} Let $(X,\delta)$ be a diversity.  Let $P_X$ denote the set of all functions
 $f \colon \Pf(X) \rightarrow \Re$ satisfying $f(\emptyset) = 0$ and
\begin{equation} \label{eqn:PXdef}
\sum_{A \in \sA} f(A) \geq \delta\big(\bigcup_{A \in \sA} A \big)
\end{equation}
for all finite $\sA \subseteq \Pf(X)$. Write $f \preceq g$ if $f(A) \leq g(A)$ for all finite $A \subseteq X$. The tight span of $(X,\delta)$ is the set $T_X$ of functions in $P_X$ that are minimal under $\preceq$. Note that if $n=|X|$ then $T_X$ can be viewed as a subset of $\Re^{(2^n-1)}$.
\end{defn}

\noindent {\bf Example. 1.} Any diversity $\delta$ on $X=\{ 1,2,3\}$ is determined by the four values
\[
d_{12}=\delta(\{1,2\}), \ \ \ \ d_{23}=\delta(\{2,3\}), \ \ \ d_{13}=\delta(\{1,3\}) , \ \ \
d_{123} = \delta(\{1,2,3\})
\]
which satisfy the condition
$$ d_{ik}\le d_{123}\le d_{ij}+d_{jk}$$
for any three distinct $i,j,k$ in $X$.

We write $f_i = f(\{i\})$, $f_{ij} = f(\{i,j\})$ and $f_{123} = f(\{1,2,3\})$ for $i,j \in X$. Condition~\eqref{eqn:PXdef} then translates to the following set of inequalities:
\begin{eqnarray}
f_i & \geq & 0 \nonumber \\
f_{ij} & \geq & d_{ij} \nonumber \\
f_i + f_j &\geq& d_{ij} \hspace{2cm} \label{eq:tight3} \\
f_{123} & \geq & d_{123} \nonumber \\
f_i + f_{jk} &\geq & d_{123}\nonumber \\
f_1 + f_2 + f_3 & \geq & d_{123} \nonumber
\end{eqnarray}
for distinct $i,j,k \in X$. Note we have omitted inequalities like $f_{ij} + f_{jk} \geq d_{123}$ since these are implied by \eqref{eq:tight3} and the triangle inequality (D2).   The elements of $T_X$ are the minimal  $f$ in $P_X$.  Equivalently, $T_X$ is the set of $f$ that satisfy 
 \eqref{eq:tight3} and such that for each nonempty $A \subseteq X$,  $f_A$ appears in an inequality in \eqref{eq:tight3} that is tight.

\begin{figure}[htbp] 
   \centering
   \includegraphics[width=4in]{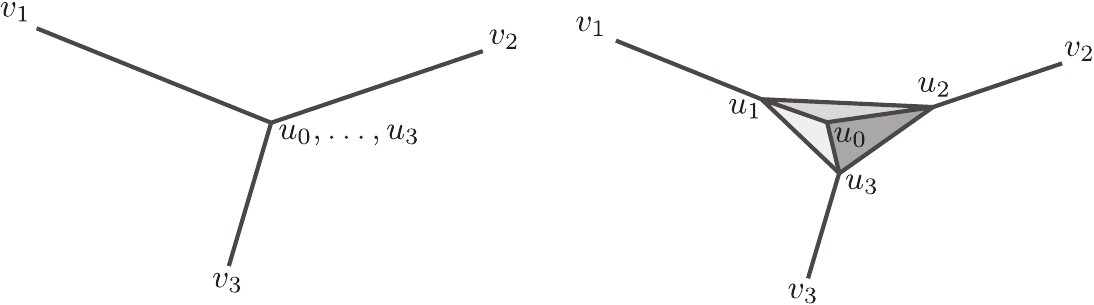} 
   \caption{Two examples of the tight span on three points, with different values for $d(\{1,2,3\})$. On the left an example where $2d_{123} \leq  d_{12} +d_{23}+d_{13}$, and the diversity tight span is one-dimensional and resembles the tight span of the induced metric. On the right a case with $2d_{123} >  d_{12} +d_{23}+d_{13}$, where the diversity consists of a 3-dimensional cell (and faces) with three adjacent 1-dimensional cells. }
   \label{fig:ThreePoint}
\end{figure}

Define the three `external' vertices
\begin{eqnarray*}
v^{(1)} &=& (0,  d_{12}, d_{13}) \\
v^{(2)} & =& (d_{12} ,  0 , d_{23} )\\
v^{(3)} &= &(d_{13},  d_{23} , 0 )
\end{eqnarray*}
and the four `internal' vertices
\begin{eqnarray*}
u^{(0)} &=& (d_{123}-d_{23},d_{123}-d_{13},d_{123}-d_{12})\\
u^{(1)} & = & u^{(0)} - (\beta,0,0) \\
u^{(2)} & = & u^{(0)} - (0,\beta,0) \\
u^{(3)} & = & u^{(0)} - (0,0,\beta),
\end{eqnarray*}
where $\beta = \max(2 d_{123} - d_{12} -d_{23}-d_{13},0)$. Let $C$ be the cell complex formed from the line segments $[u^{(1)},v^{(1)}]$, $[u^{(2)},v^{(2)}]$, $[u^{(3)},v^{(3)}]$ and the solid tetrahedron with vertices $u^{(1)},\ldots,u^{(4)}$.  We will show that  $f \in T_X$ if and only if $(f_1,f_2,f_3) \in C$, 
$f_{23}  =  \max(  d_{23}, d_{123} - f_1 )$, $f_{13}   =  \max( d_{13}, d_{123} - f_2 )$,  $f_{23}  =  \max( d_{12}, d_{123} - f_3)$, $f_{123}=d_{123}$. If $\beta = 0$ then $u^{(0)}$ to $u^{(3)}$ coincide, and the tight span is one-dimensional and resembles the metric tight span for the induced metric, albeit sitting in $\Re^7$ (Figure~\ref{fig:ThreePoint}a).
When $\beta > 0$ the tight span resembles a tetrahedron with three spindles branching off, as in Figure~\ref{fig:ThreePoint}b.

To show that this is the tight span,
first note that $f_{12} = \max(d_{12},d_{123}-f_3)$, so $f_{12}$ is fixed once $f_1,f_2,f_3$ are fixed. By symmetry, the same holds of $f_{13}$ and $f_{23}$. Moreover, $f_{123}=d_{123}$.  So points in the tight span are uniquely characterized by the values of $f_1,f_2,f_3$. 
\begin{itemize}
\item 
\noindent {\em Case 1.} At least two of the inequalities $f_i+f_j \geq d_{ij}$ are tight.\\ 
  In this case $(f_1,f_2,f_3)$ is in the metric tight span, with the additional constraint that $f_1+f_2+f_3 \geq d_{123}$.  These values of $f_1,f_2,f_3$ correspond to the line segments $[u^{(1)},v^{(1)}]$, $[u^{(2)},v^{(2)}]$, $[u^{(3)},v^{(3)}]$.
\item
\noindent{\em Case 2}. At most one of the inequalities $f_i+f_j \geq d_{ij}$ is tight. \\
Without loss of generality, suppose that $f_1+f_2>d_{12}$ and $f_1+f_3>d_{13}$. If $f_1>d_{123} - d_{23}$ then $f_1 + f_{23} > d_{123}$ and since $f_2 + f_3 \geq d_{23}$ we would also obtain $f_1+f_2+f_3 > d_{123}$. This then leaves no tight inequality involving $f_1$. Hence we conclude $f_1 \leq d_{123}-d_{23}$.

A similar analysis shows that we also have $f_2 \leq d_{123} - d_{13}$ and $f_3 \leq d_{123} - d_{12}$
as assuming for example that $f_2 >d_{123} - d_{13}$ holds, we would also have Ê$f_2+f_{13} \ge Êf_2+d_{13}> Êd_{123}$ as well as
$f_1+f_2+f_3 \ge f_2+d_{13} > d_{123}$ and hence, in view of Ê$f_1 \leq d_{123}-d_{23}$, also
$ f_2+f_3 > Êd_{123}-f_1\ge d_{23}$ leaving no tight inequality involving $f_2$."

We have now shown that $f_i \leq d_{123} - d_{jk}$ for distinct $i,j,k$. These inequalities together with $f_1+f_2+f_3 \geq d_{123}$ define the tetrahedron given by the vertices $u^{(0)},u^{(1)},u^{(2)},u^{(3)}$. \hfill $\square$\\
\end{itemize}

We now prove a characterisation of the diversity tight span which will be used extensively throughout the remainder of the paper  (Theorem~\ref{thm:altchartight}). An equivalent result holds for the metric tight span \cite[Theorem 3(v)]{Dress84}. 

\begin{thm} \label{thm:altchartight}
Let $f \colon \Pf(X) \rightarrow \Re$ and suppose $f(\emptyset) = 0$. Then $f \in T_X$ if and only if for all finite $A \subseteq X$,
\begin{equation}
f(A) = \sup_{\sB \subseteq \Pf(X)} \left\{\delta\big(A \cup \bigcup_{B \in \sB}B \big) - \sum_{B \in \sB} f(B) \colon |\sB| < \infty \right\}. \label{eq:altchartight}
\end{equation}
\end{thm}

\begin{pf} 
Suppose that $f \in T_X$.  
For all finite $A\subseteq X$ and all finite  $\sB \subseteq \Pf(X)$ we have 
\[f(A) \geq \delta(A \cup \bigcup_{B \in \sB} B)-\sum_{B \in \sB}f(B),\]
 giving the required lower bound on $f(A)$.  
Now suppose that for some finite $A_0$ 
\begin{equation}
f(A_0) >  \sup_{\sB \subseteq \Pf(X)} \left\{\delta\big(A_0 \cup \bigcup_{B \in \sB}B\big) - \sum_{B \in \sB} f(B) \colon |\sB| < \infty \right\}. \label{eq:BadA0}
\end{equation}
Define a function $g:\Pf(X) \rightarrow \Re_{\geq 0}$ by
\[g(A) = \begin{cases} 
f(A) & \mbox{ if $A \neq A_0$}\\
 \sup_{\sB \subseteq \Pf(X)} \left\{\delta\big(A_0 \cup \bigcup_{B \in \sB}B\big) - \sum_{B \in \sB} f(B) \colon |\sB| < \infty \right\}& \mbox{ if $A = A_0$.}
\end{cases}\]
Clearly $g \neq f$ and $g \preceq f$. We show that $g$ is in $P_X$. Let $\sA$ be a finite subset of $\Pf(X)$. If $A_0 \not \in \sA$ then
\[\sum_{A \in \sA} g(A) = \sum_{A \in \sA} f(A) \geq \delta\big(\bigcup_{A \in \sA} A \big).\]
If $A_0 \in \sA$ then
\begin{eqnarray*}
\sum_{A \in \sA} g(A) &=&  \sup_{\sB \subseteq \Pf(X)} \left\{\delta\big(A_0 \cup \bigcup_{B \in \sB}B\big) - \sum_{B \in \sB} f(B) \colon |\sB| < \infty \right\} + \sum_{B \in \sA \setminus \{A_0\}} f(B)  \\
& \geq &  \delta\big(A_0 \cup \bigcup_{B \in \sA \setminus \{A_0\}}B\big)\\
&=&  \delta( \bigcup_{A \in \sA} A) ,
\end{eqnarray*}
by letting $\sB = \sA \setminus \{ A_0 \}$.
 So $g \in P_X$, $g \neq f$ and  $g \preceq f$, contradicting $f \in T_X$. Hence there is no $A_0$ satisfying~\eqref{eq:BadA0}. If $f \in T_X$ then \eqref{eq:altchartight} holds for all finite $A \subseteq X$.

For the converse, suppose that \eqref{eq:altchartight} holds for all  finite $A \subseteq X$. Then $f \in P_X$. Suppose that $g \in P_X$, that $g \preceq f$ and $A \in \Pf(X)$. Then for all finite $\sB \subseteq \Pf(X)$ 
we have
\begin{eqnarray*}
\delta\big(A \cup \bigcup_{B \in \sB}B \big) - \sum_{B \in \sB} f(B) &\leq& \delta\big(A \cup \bigcup_{B \in \sB}B \big) - \sum_{B \in \sB} g(B)   \leq g(A)
\end{eqnarray*}
so that $f(A)  \leq g(A)$. Hence $f$ is minimal in $P_X$.
\end{pf}

We note that the characterisation of tight spans given by Theorem~\ref{thm:altchartight} is analogous to the definition of tight spans for valuated matroids used by \cite{Dress98}. One important difference is that, for diversities, the tight span is made up of functions on $\Pf(X)$ rather than functions on $X$.

The following basic properties of members of $T_X$ will be used subsequently. 

\pagebreak

\begin{prop} \label{cor:tightmonotone}
Suppose that $f \in T_X$. 
\begin{enumerate}
\item $f(A) \geq \delta(A)$ for all finite $A \subseteq X$.
\item If $A \subseteq B \subseteq X$ and $B$ is finite then $f(A) \leq f(B)$; that is, $f$ is {\em monotone}.
\item $f(A \cup C) \leq \delta(A \cup B) + f(B \cup C)$ for all $A,B,C \in \Pf(X)$ with $B \neq \emptyset$.
\item $f(A \cup B) \leq f(A) + f(B)$ for all $A,B \in \Pf(X)$; that is, $f$ is {\em sub-additive}.
\item For all finite $A$,
\begin{equation}
f(A) = \sup_B\left\{ \delta( A \cup B) - f(B) : B \in \Pf(X) \right\}. \label{eq:TwoWayTight}
\end{equation}
\end{enumerate}
\end{prop}
\begin{pf}
1. Use $\sA = \{A\}$ in the definition of $P_X$.\\
2. Follows from \eqref{eq:altchartight} and the monotonicity of $\delta$.  \\
%
3. Let $A,B, C \in \Pf(X)$ with $B \neq \emptyset$.
We have
\begin{eqnarray*}
f(A \cup C) & = Ê& \hspace*{-0.3cm}\sup_{\sD \subseteq \Pf(X) }\left\{\delta(A \cup C \cup \bigcup_{D \in \sD} D) Ê- \sum_{D \in \sD} f(D) \colon | \sD| < \infty \right\} Ê\\
& \le & \hspace*{-0.3cm}\sup_{\sD \subseteq \Pf(X) }\left\{\delta(A \cup B)+ \delta(B \cup C \cup \bigcup_{D \in \sD} D) Ê- \sum_{D \in \sD} f(D) \colon | \sD| < \infty \right\}\\
Ê& = Ê& \hspace*{-0.3cm}\delta(A \cup B)+\sup_{\sD \subseteq \Pf(X) }\left\{ \delta(B \cup C\cup Ê\bigcup_{D \in \sD} D) Ê- \sum_{D \in \sD} f(D) \colon | \sD| < \infty \right\}\\
ÊÊ& = Ê& \hspace*{-0.3cm}\delta(A \cup B)+f(B \cup C )\\
\end{eqnarray*}
We note that this property is analogous to the continuity of functions in the metric tight span, see \cite[Theorem 3(iv)]{Dress84}. \\
4. Given any $A,B \in \Pf(X)$ and any finite collection $\sC \subseteq \Pf(X)$ we have
\[f(A) + f(B) + \sum_{C \in \sC} f(C) \geq \delta\big(A \cup B \cup \bigcup_{C \in \sC}C \big) \]
so that
\begin{eqnarray*}
f(A) + f(B) & \geq & \sup_{\sC \subseteq \Pf(X)} \left\{ \delta\big(A \cup B \cup \bigcup_{C \in \sC}C \big) -\sum_{C \in \sC} f(C)\colon | \sC| < \infty \right\} \\
& = & f(A \cup B)
\end{eqnarray*}
by Theorem~\ref{thm:altchartight}. \\
5. For any finite $\sB \subseteq \Pf$,
\(
\sum_{B \in \sB} f(A) \geq f( \cup_{B \in \sB} B).
\)
So
\[
\sup_{\sB \in \Pf(X) } \left\{ \delta( A \cup \bigcup_{B \in \sB} B ) - \sum_{B\in \sB} f(B) \colon |\sB| < \infty \right\} = \sup_{C \in \Pf(X) } \left\{ \delta( A \cup C) - f(C) \right\}.
\]
\end{pf}
It was recently shown in \cite{Herrman11} that if $(X,\delta)$ is an $L_1$ diversity then \eqref{eq:TwoWayTight} is sufficient for $f \in T_X$. This does not hold in general, even when $|X|=3$.\\

\noindent {\bf Example 2.} Consider $X = \{1,2,3\}$ and the diversity $\delta$ defined as in Example 1.~with $d_{12} = d_{13} = d_{23} = 1$ and $d_{123} = 2$. Define $f$ (using the same notation) with $f_1 = f_2 = f_3 = 1/2$, $f_{12}=f_{13}=f_{23} = 3/2$ and $f_{123} = 2$. Then $f$ satisfies \eqref{eq:TwoWayTight} but it is not in $T_X$, since $f_1+f_2+f_3<d_{123}$.  \hfill $\square$\\

The distance between any two functions $f,g$ in the metric tight span is given by the $l_\infty$ norm, 
\begin{equation}
d_T(f,g) = \sup_{x \in X} |f(x) - g(x)| \label{eq:linfty}
\end{equation}
which Dress \cite[Theorem 3(iii)]{Dress84} shows is equivalent on this set to 
\begin{equation} \label{eq:dTdef}
d_T(f,g) = \sup_{x,y \in X} \{d(x,y) - f(x) - g(y)\}.
\end{equation}
Dress  also showed that a metric can be embedded into its tight span using the Kuratowski map $\kappa$, which takes an element $x \in X$ to the function  $h_x$ for which $h_x(y)=d(x,y)$ for all $y$. This is exactly the map $e$ defined in \cite[section 2.4]{isbell}.

Here we establish the analogous results for the diversity tight span. We define the appropriate function $\kappa$ from a diversity to its tight span. We then define a  function $\delta_T$ on $T_X$ so that $(T_X,\delta_T)$ is a diversity and prove that $\kappa$ is an embedding. 

\begin{defn}
 \begin{enumerate} 
\item Let $(Y_1,\delta_1)$ and $(Y_2,\delta_2)$ be two diversities.  A map $\pi \colon Y_1 \rightarrow Y_2$ is an embedding if it is  one-to-one (injective) and for all finite $A \subseteq Y_1$ we have $\delta_1(A) = \delta_2(\pi(A))$.  In this case, we say that $\pi$ embeds $(Y_1,\delta_1)$ in $(Y_2,\delta_2)$.
\item An isomorphism is an onto (surjective) embedding between two diversities.
\item 
Let $(X,\delta)$ be a diversity. For each $x \in X$ define the function \(h_x:\Pf(X) \rightarrow \Re \) by \[ h_x(A) = \delta(A \cup \{x\})\] for all finite $A \subseteq X$. Let $\kappa$ be the map taking each $x \in X$ to the corresponding function $h_x$.
\item
Let $(X,\delta)$ be a diversity. 
 Let $\delta_T:\Pf(T_X) \rightarrow \Re$  be the function defined by $\delta_T(\emptyset) = 0$ and 
\begin{eqnarray}
\delta_T(F)  &=& \sup_{\sA \subseteq \Pf(X)} \left\{ \delta\left( \bigcup_{A \in \sA} A \right) - \sum_{A \in \sA} \inf_{f \in F} f(A) : |\sA|<\infty \right\} \label{eq:deltaTdef}  
\end{eqnarray}
for all finite non-empty $F \subseteq T_X$. 
\end{enumerate}
 \end{defn}


Further manipulations give a form for $\delta_T$ analogous to \eqref{eq:dTdef}:
\[
\delta_T(F) = \sup_{ \{A_f\}_{f \in F}} \left\{ \delta \left( \bigcup_{f \in F} A_f \right)  - \sum_{f \in F} f(A_f) : A_f \in \Pf(X) \mbox{ for all $f \in F$}\right\},
\]
for all finite $F \subseteq \Pf(T_X)$.
We can also re-express \eqref{eq:deltaTdef} in a form closer to \eqref{eq:linfty}. Note the similarity between Lemma~\ref{lem:deltat} and \cite[Theorem 3(iii)]{Dress84}.

\begin{lem} \label{lem:deltat} 
If $f \in F$ then
\begin{eqnarray*}
\delta_T(F) & = & \sup_{\sA \subseteq \Pf(X)} \left\{ f\left( \bigcup_{A \in \sA}A\right) - \sum_{A \in \sA} \inf_{g \in F \setminus \{f\}} g(A) : |\sA|<\infty\right\} .
\end{eqnarray*}
\end{lem}
\begin{pf}
For $\sA \subseteq \Pf(X)$ define
 \begin{align*} \sA' &= \{A \in \sA : f(A) > \inf_{g \in F} g(A)\}\\
 \sA''& = \{A \in \sA : f(A) = \inf_{g \in F} g(A)\}. \end{align*}
Then
\begin{align*}
\delta_T(F) & = \sup_{\sA \subseteq \Pf(X)} \left\{ \delta\left( \bigcup_{A \in \sA} A \right) - \sum_{A \in \sA} \inf_{g \in F} g(A) : |\sA|<\infty \right\} \\
& = \sup_{\sA \subseteq \Pf(X)} \left\{ \delta\left( \bigcup_{A \in \sA} A \right) - \sum_{A \in \sA'} \inf_{g \in F\setminus \{f\}} g(A) - \sum_{A \in \sA''} f(A) : |\sA|<\infty \right\}\\
& = \sup_{\sB,\sC \subseteq \Pf(X)} \left\{ \delta\left(\left( \bigcup_{B \in \sB} B \right) \cup \left( \bigcup_{C \in \sC} C \right)   \right) - \sum_{B \in \sB} \inf_{g \in F\setminus \{f\}} g(B) - \sum_{C \in \sC} f(C) : |\sB|,|\sC| <\infty \right\}. \\
\intertext{This last line follows from that fact that if $B \in \sB$ and $f(B) = \inf_{g \in F} g(B)$ then moving $B$ from $\sB$ to $\sC$ cannot decrease 
\begin{equation}
\delta\left(\left( \bigcup_{B \in \sB} B \right) \cup \left( \bigcup_{C \in \sC} C \right)   \right) - \sum_{B \in \sB} \inf_{g \in F\setminus \{f\}} g(B) - \sum_{C \in \sC} f(C) \label{eq:otherdt}
\end{equation}
while if $C \in \sC$ and $f(C) > \inf_{g \in F} g(C)$ then moving $C$ from $\sC$ to $\sB$ will increase \eqref{eq:otherdt}. Continuing, we have }
\delta_T(F) & =  \sup_{\substack{\sB \subseteq \Pf(X)\\ |\sB|<\infty}} \left\{ \sup_{\substack{\sC \subseteq \Pf(X)\\|\sC|<\infty}} \left\{ \delta\left(\left( \bigcup_{B \in \sB} B \right) \cup \left( \bigcup_{C \in \sC} C \right)   \right) - \sum_{C \in \sC} f(C) \right\} - \sum_{B \in \sB} \inf_{g \in F\setminus \{f\}} g(B) \right\}. \\
& = \sup_{\sB \subseteq \Pf(X)} \left\{f\left( \bigcup_{B \in \sB} B\right) - \sum_{B \in \sB} \inf_{g \in F \setminus \{f\}} g(B)  : |\sB|< \infty  \right\}.
\end{align*}
by Theorem~\ref{thm:altchartight}.\end{pf}


\begin{thm}  \label{thm:tightspan_is_diversity}
$(T_X,\delta_T)$ is a diversity. 
\end{thm}
\begin{pf}
First note that for all $F \subseteq T_X$, when $\sA = \{\emptyset\}$,
\[\deltap{ \bigcup_{A \in \sA}A}  - \sum_{A \in \sA} \inf_{f \in F} f(A) = 0\]
so that $\delta_T$ is non-negative.

If $\emptyset \neq F \subseteq G$ then for all $\sA \subseteq \Pf(X)$ with $|\sA|<\infty$ we have \[\sum_{A \in\sA} \inf_{f \in F} f(A) \geq \sum_{A \in \sA} \inf_{f \in G} f(A).\]
Hence $\delta_T(F) \leq \delta_T(G)$, showing that $\delta_T$ is monotone.

If $F = \{f\}$ then
\[\delta_T(F) = \delta_T(\{f\}) \leq \sup \left\{\deltap{ A} - f(A) :A \in \Pf(X) \right\} =0\]
 by the subadditivity of $f$ and by part 1 of Proposition~\ref{cor:tightmonotone}. 
On the other hand, if $|F|>1$ then there is $f_1,f_2 \in F$ such that $f_1 \neq f_2$.  By monotonicity and 
  Lemma~\ref{lem:deltat} we have
\begin{eqnarray*}
\delta_T(F) & \geq & \delta_T(\{f_1,f_2\}) \\
& = & \sup_{\sA \subseteq \Pf(X)} \left\{ f_1\left(\bigcup_{A \in \sA} A \right) - f_2\left(\bigcup_{A \in \sA} A \right) : |\sA|<\infty \right\} \\
& = & \sup_{A \in \Pf(X)} \left\{ f_1(A) - f_2(A) \right\} \\
& >  0.&
\end{eqnarray*}
We have now proved that $\delta_T$ satisfies (D1).

For the triangle inequality, suppose $F$ and $G$ are disjoint finite subsets of $T_X$ and that $h \in T_X \setminus (F \cup G)$.  Then by Lemma~\ref{lem:deltat}
\begin{eqnarray} \label{eq:dt1}
\delta_T(F \cup \{h\}) & = & \sup_{\sA \subseteq \Pf(X)} \left\{h \left( \bigcup_{A \in \sA}A  \right) - \sum_{A \in \sA} \inf_{f \in F} f(A) : |\sA|<\infty \right\}
\end{eqnarray}
and
\begin{eqnarray}
\delta_T(G \cup \{h\}) & = & \sup_{\sB \subseteq \Pf(X)} \left\{h \left( \bigcup_{B \in \sB}B \right) - \sum_{B \in \sB} \inf_{g \in G} g(B) : |\sB|<\infty \right\}. \label{eq:dt2}
\end{eqnarray}
By part 4 of Proposition~\ref{cor:tightmonotone} the function $h$ is sub-additive, so 
\begin{equation}
h \left( \bigcup_{A \in \sA}A  \right) + h \left( \bigcup_{B \in \sB} B \right) 
\geq h \left(   \bigcup_{C \in \sA \cup \sB }C  \right) .\label{eq:dt3}
\end{equation}
Combining \eqref{eq:dt1}--\eqref{eq:dt3} and again applying Lemma~\ref{lem:deltat} we have
\begin{eqnarray*}
\delta_T(F \cup \{h\}) + \delta_T(G \cup \{h\}) & \geq &
\sup_{\sC \subseteq \Pf(X)} \left\{ h\left( \bigcup_{C \in \sC }C \right) - \sum_{C \in \sC} \inf_{f \in F \cup G} f(C) :|C|<\infty \right\} \\
& = & \delta_T(F \cup G \cup \{h\}).
\end{eqnarray*}
The triangle inequality (D2) now follows by monotonicity.

\end{pf}

Theorem~\ref{thm:tightspan_is_diversity} establishes that $(T_X,\delta_T)$ is a diversity. We now show that $\kappa$ is an embedding from $(X,\delta)$ into $(T_X,\delta_T)$. We then prove the diversity analogue of \cite[Eq. (2.4)]{isbell} (see \cite[Theorem 3(ii)]{Dress84}) and  characterise $\delta_T$ in terms of a minimality condition. 

\begin{thm} \label{thm:deltaTmin}
\begin{enumerate}
\item The map $\kappa$ is an embedding from $(X,\delta)$ into $(T_X,\delta_T)$. 
\item For all finite $Y \subseteq X$  and $f \in T_X$, 
\[\delta_T(\kappa(Y) \cup \{f\}) = f(Y).\]
\item If $(T_X,\hdelta)$ is a diversity such that $\hdelta(\kappa(Y) \cup \{f\}) = f(Y)$ for all finite $Y \subseteq X$ and $f \in T_X$ then
\[\hdelta(F) \geq \delta_T(F)\]
for all finite $F \subseteq T_X$.
\end{enumerate}
\end{thm}
\begin{pf}
1. 
 Fix $x \in X$. Consider finite $\sA \subseteq \Pf(X)$. The triangle inequality for diversities, (D2), gives
\[\sum_{A \in \sA} h_x(A) = \sum_{A \in \sA} \delta(A \cup \{x\}) \geq \delta\left(\bigcup_{A \in \sA} A \right),\]
so that $h_x \in P_X$. There is $g \in T_X$ such that $g \preceq h_x$. Since $h_x(\{x\}) = \delta(\{x\}) = 0$ we have for all finite $A \subseteq X$ that
\[h_x(A)\: = \: \delta(A \cup \{x\}) \:\:\: \leq \:\:\: g(A) \!+\! g(\{x\}) \:\:\: \leq \:\:\: g(A) \!+\! h_x(\{x\}) \:=\: g(A) \:\: \leq \:\:\: h_x(A).\]
Hence $h_x = g \in T_X$.

To see that $\kappa$ is one-to-one observe that for $x \neq y$, $h_x(\{x\}) = 0$ but $h_y(\{x\})=\delta(\{x,y\}) >0$.  So $h_x \neq h_y$ for distinct $x,y \in X$.

We now show that $\delta_T(\kappa(Y)) = \delta(Y)$ for all finite $Y \subseteq X$. Let $Y \subseteq X$, $Y=\{y_1, \ldots, y_k\}$.  Taking $\sA = \{ \{y_1\},\ldots, \{y_k\}\}$ in  \eqref{eq:deltaTdef} gives $\delta_T(\kappa(Y)) \geq \delta(Y)$. 

By repeatedly using the triangle inequality we have for any finite $\sA = \{ A_1, A_2, \ldots, A_j \} \subseteq \Pf(X)$ and $z_1, \ldots, z_j \in Y$ that 
\begin{eqnarray*}
\delta(Y) & \geq & \delta( Y  \cup A_1 ) - \delta( \{z_1\} \cup A_1) \\
 & \geq & \delta( Y \cup  A_1 \cup A_2 ) - \delta( \{z_1\} \cup A_1) - \delta( \{z_2 \} \cup A_2) \\
 & \geq & \delta\left( Y \cup \bigcup_{i=1}^j A_i \right) - \sum_{i=1}^j \delta( \{z_i\} \cup A_i) \\
 & \geq & \delta\left( \bigcup_{i=1}^j A_i \right) - \sum_{i=1}^j h_{z_i}(  A_i) \\
 & \geq & \delta\left( \bigcup_{i=1}^j A_i \right) - \sum_{i=1}^j \inf_{h \in \kappa(Y)}h(A_i) .
 \end{eqnarray*}
 Taking the supremum over all such $\sA$  and applying \eqref{eq:deltaTdef} gives
 $\delta_T(\kappa(Y)) \leq \delta(Y)$.
So $\delta_T(\kappa(Y))= \delta(Y)$ and $\kappa$ is an embedding.

2. Let $Y \subseteq X$, $Y$ finite, and $f \in T_X$. If $f=h_y$ for $y\in Y$ then, using part 1, 
\[
\delta_T( \kappa(Y) \cup \{f\} ) = \delta_T(\kappa(Y)) = \delta( Y) = \delta( Y \cup \{y\}) = f(Y),
\]
as required.  Otherwise,  suppose $f \not\in \kappa(Y)$.
 Let $Y = \{ y_1,\ldots, y_k \}$, so that
\[
\delta_T( \kappa(Y) \cup \{f\} ) = \sup_{A_i, i=1,\ldots, k, A_f} \left\{ \delta\left( \bigcup_{i} A_i  \cup A_f\right) - \sum_i \delta(\{y_i\} \cup A_i)  - f(A_f) \right\}.
\]
Letting $A_i = \{y_i\}$ for all $i$ shows 
\[
\delta_T(\kappa(Y) \cup \{f\} ) \geq \sup_{A_f} \{\delta( Y \cup A_f ) - f(A_f)\} = f(Y),
\]
by Proposition~\ref{cor:tightmonotone} part 5. On the other hand, following the same reasoning as in part 1 of this proof shows
\[
\delta_T(\kappa(Y) \cup \{f\})  \leq \sup_{A_f} \delta(Y \cup A_f ) - f(A_f) = f(Y).
\]
Therefore $\delta_T(\kappa(Y) \cup \{f\}) = f(Y)$. 

3. Suppose that $F = \kappa(Y) \cup G$, where $Y \in \Pf(X)$ and $G \subseteq T_X \setminus \kappa(X)$. For all collections $\sA \subseteq \Pf(X)$ with $|\sA|<\infty$ and all collections $\{f_A\}_{A \in \sA}$ of elements in $F$, we have from parts 1 and 2 that
\begin{eqnarray*}
\deltap{Y \cup \bigcup_{A \in \sA} A} - \sum_{A \in \sA} f_A(A) & = & \hdelta\left( \kappa(Y) \cup \bigcup_{A \in \sA} \kappa(A) \right) - \sum_{A \in \sA} \hdelta(\kappa(A) \cup \{f_A\}) \\ 
& \leq & \hdelta\left(\kappa(Y) \cup \{f_A:A \in \sA\} \right) \\
& \leq &  \hdelta\left(\kappa(Y) \cup F \right).
\end{eqnarray*}
\end{pf}

\section{Hyperconvex diversities and the injective envelope} \label{sec:hyperconvex} 

Aronszajn and Panitchpakdi \cite{Aronszajn56} introduced hyperconvex metric spaces and showed that they are exactly the injective metric spaces. 
\begin{defn} \label{def:metricTS}
\begin{enumerate}
\item A metric space $(X,d)$ is said to be {\em hyperconvex} if for all $r \colon X \rightarrow \Re$ with $r(x) + r(y) \geq d(x,y)$ for all $x,y \in X$ there is a point $z \in X$ such that $d(z,x) \leq r(x)$ for all $x \in X$.
\item A metric space $(X,d)$ is {\em injective} if it satisfies the following property: given any pair of metric spaces $(Y_1,d_1)$, $(Y_2, d_2)$, an embedding $\pi:Y_1 \rightarrow Y_2$ and a non-expansive map $\phi:Y_1 \rightarrow X$ there is a non-expansive map 
$\psi:Y_2 \rightarrow X$ such that $\phi = \psi \circ \pi$.
\end{enumerate}
\end{defn}
See \cite{Espinola01} for a proof of the equivalence of these two concepts, as well as a highly readable and comprehensive review of the rich metric structure of hyperconvex spaces. Here we establish diversity analogues for these concepts and show that the equivalence holds in this new setting. 
We begin by defining diversity analogues of injective and hyperconvex  metric spaces.

\begin{defn}
\begin{enumerate}
\item Given diversities $(Y_1,\delta_1)$ and $(Y_2,\delta_2)$, a map $\phi:Y_1 \rightarrow Y_2$ is {\em non-expansive} if for all $A \subseteq Y_1$ we have $\delta_1(A) \geq \delta_2(\phi(A))$ and it is an {\em embedding} if it is one-to-one and for all $A \subseteq Y_1$ we have $\delta_1(A) =\delta_2(\phi(A))$. 
\item  A diversity $(X,\delta)$ is {\em injective} if it satisfies the following property: given any pair of diversities $(Y_1,\delta_1)$, $(Y_2, \delta_2)$, an embedding $\pi:Y_1 \rightarrow Y_2$ and a non-expansive map $\phi:Y_1 \rightarrow X$ there is a non-expansive map 
$\psi:Y_2 \rightarrow X$ such that $\phi = \psi \circ \pi$.
\item 
A diversity $(X,\delta)$ is said to be {\em hyperconvex} if for all $r \colon \Pf(X) \rightarrow \Re$ such that 
\begin{equation} \label{hypercon_cond}
\delta \left( \bigcup_{A \in \sA} A \right) \leq \sum_{A \in \sA} r(A) 
\end{equation}
for all finite $\sA \subseteq \Pf(X)$ there is $z \in X$ such that
$\delta( \{ z\} \cup Y) \leq r(Y)$ for all finite $Y \subseteq X$.
\end{enumerate}
\end{defn}

The following theorem establishes the diversity equivalent of Aronszajn and Panitchpakdi's result.  

\begin{thm} \label{lem:injecthyper}
A diversity $(X,\delta)$ is injective if and only if it is hyperconvex.
\end{thm}
\begin{pf}
First suppose that $(X,\delta)$ is injective. Consider $r \colon \Pf(X) \rightarrow \Re$ satisfying \eqref{hypercon_cond} for all finite $\sA \subseteq \Pf(X)$.  Without loss of generality we can assume $r(\emptyset)=0$ and hence $r \in P_X$.   Choose $f \in T_X$ with $f \preceq r$.

Let $x^*$ be a point not in $X$, let $X^* = X \cup \{x^*\}$  and let $\delta^*:\Pf(X \cup \{x^*\}) \rightarrow \Re$ be the function where for all finite $A \subseteq X$,
\begin{eqnarray*}
\delta^*(A) & = & \delta(A) \\
\delta^*(A \cup \{x^*\}) & = & f(A).
\end{eqnarray*}
From part 2 of Proposition~\ref{cor:tightmonotone} we have that $\delta^*$ is monotonic, and from parts 4 and 5 we have that
\begin{eqnarray}
\delta^*(A \cup C \cup \{x^*\}) &\leq & \delta^*(A \cup \{x^*\}) + \delta^*(C \cup \{x^*\}) \label{eq:injHyper1} \\
\delta^*(A \cup B \cup C \cup \{x^*\}) &\leq & \delta^*(A \cup B \cup \{x^*\}) + \delta^*(B \cup C) \label{eq:inHyper2} 
\end{eqnarray}
for all finite $A,B,C \subseteq X$ such that $B \neq \emptyset$. These, together with monotonicity and the fact that $\delta^*$ coincides with $\delta$ on $\Pf(X)$, imply the triangle inequality (D2) for $(X^*,\delta^*)$. 

We now apply the fact that $(X,\delta)$ is injective. Let $(Y_1,\delta_1)$ be $(X,\delta)$, let $(Y_2,\delta_2)$ be $(X^*,\delta^*)$, let $\pi$ be the identity embedding from $(X,\delta)$ into $(X^*,\delta^*)$ and let $\phi$ be the identity map from $(X,\delta)$ to itself. Then there is a non-expansive map $\phi:X^* \rightarrow X$ such that $\phi(x) = x$ for all $x \in X$.

Let $\omega = \phi(x^*)$. For all finite $A \subseteq X$ we have
\begin{eqnarray*}
\delta(A \cup \{\omega\}) & \leq & \delta^*(A \cup \{x^*\}) \\
& = & f(A) \\
& \leq & r(A).
\end{eqnarray*}
This proves that $(X,\delta)$ is hyperconvex.

For the converse, suppose now that $(X,\delta)$ is hyperconvex. Let $(Y_1,\delta_1)$ and $(Y_2,\delta_2)$ be two diversities, let $\pi:Y_1 \rightarrow Y_2$ be an embedding and let $\phi$ be a non-expansive map from $Y_1$ to $X$.  We will show that there is non-expansive $\psi:Y_2 \rightarrow X$ such that $\phi = \psi \circ \pi$.

Let $\sY$ denote the collection of pairs $(Y,\psi_Y)$ such that  $\pi(Y_1) \subseteq Y \subseteq Y_2$  and  $\psi_Y$ is a non-expansive map from $Y$ to $X$  such that $\phi = \psi_Y \circ \pi$.   We want to show that $Y_2 \in \sY$. Suppose this is not the case.  We write $(Y,\psi_Y) \unlhd (Z,\psi_Z)$ if $Y \subseteq Z$ and $\psi_Z$ restricted to $Y$ equals $\psi_Y$. The partially ordered set $(\sY, \unlhd)$ satisfies the conditions of Zorn's lemma, so it contains maximal elements. 

Let $(Y,\psi_Y)$ be one such maximal element. Choose $y \in Y_2 \setminus Y$. For each finite 
$A \subseteq Y$ let $r(A) = \delta_2(A \cup \{y\})$. For any finite collection $\sA \subseteq \Pf(Y)$ 
we have
\begin{eqnarray*}
\delta\left(\bigcup_{A \in \sA} \psi_Y(A)\right) & = & \delta\left(\psi_Y \left( \bigcup_{A \in \sA} A \right) \right) \\ 
& \leq & \delta_2\left( \bigcup_{A \in \sA} A \right)  \\  
& \leq & \sum_{A \in \sA}  \delta_2(A \cup \{y\}) \\
& = & \sum_{A \in \sA} r(A).
\end{eqnarray*}
If $A \not \subseteq \psi_Y(Y)$ we let $r(A) = \infty$. Since $(X,\delta)$ is hyperconvex, there is $x\in X$ such that 
\[\delta(\psi_Y(A) \cup \{x\}) \leq r(A) = \delta_2(A \cup \{y\})\]
for all finite  $A \subseteq Y$. Hence we can extend $\psi_Y$ to $Y \cup \{y\}$ by setting $\psi_Y(y) = x$, giving a non-expansive map from $Y \cup \{y\}$ to $X$, and contradicting the maximality of $Y$.

It follows that $Y_2 \in \sY$, proving that $(X,\delta)$ is injective.
\end{pf}

\begin{defn}
Let $(X,\delta)$ be a diversity. For $F \subseteq T_X$ and finite $Y \subseteq X$ let 
\[\Phi_F(Y) = \inf_{\sA \subseteq \Pf(X)} \left\{ \sum_{A \in \sA} \inf_{f \in F} f(A ) : |\sA| < \infty,
\bigcup_{A \in \sA} A = Y
 \right\}.\]
\end{defn}

Clearly, 
\begin{equation}
\delta_T(F)  = \sup_{Y \subseteq X} \{\delta(Y) - \Phi_F(Y) \colon |Y|<\infty \}.\label{eq:deltatPhi}
\end{equation}
We show that $\Phi_F$ also satisfies a sub-additivity type identity.

\begin{lem} \label{lem:phi}
For $F,G \subseteq T_X$ and $Y,Z \subseteq \Pf(X)$ we have
\[\Phi_{F \cup G} (Y \cup Z) \leq \Phi_F(Y) + \Phi_G(Z) .\]
\end{lem}
\begin{pf}
Given $\epsilon > 0$ there is finite $\sA \subseteq \Pf(X)$  and a collection $\{f_A\}_{A \in \sA}$ of elements in $T_X$ such that
\[\Phi_{F}(Y) \,\, \leq \,\, \sum_{A \in \sA} f_A(A ) \,\, < \,\, \Phi_{F}(Y) + \epsilon/2.\]
Similarly, there is finite $\sB \subseteq \Pf(X)$  and a collection $\{g_B\}_{B \in \sB}$ of elements in $T_X$ such that
\[\Phi_{G}(Z) \,\,\leq \,\, \sum_{B \in \sB} g_B(B ) \,\, <\,\,  \Phi_{G}(Z) + \epsilon/2.\]
Define $\sC = \sA \cup \sB$ and the collection $\{h_C\}_{C \in \sC}$ by
\[h_C = \begin{cases} f_C & \mbox{ if $C \in \sA$;} \\ g_C & \mbox{ otherwise. } \end{cases}\]
Then 
\begin{eqnarray*}
\Phi_F(Y) + \Phi_G(Z) + \epsilon & > &  \sum_{A \in \sA} f_A (A ) + \sum_{B \in \sB} g_B (B ) \\
& \geq & \sum_{C \in \sC} h_C(C) \\
& \geq & \Phi_{F \cup G} (Y \cup Z).
\end{eqnarray*}
Taking $\epsilon \rightarrow 0$ proves the lemma.
\end{pf}

Isbell  proved that the metric tight span is injective, and hence hyperconvex \cite[Section 2.9]{isbell}. Here we prove the same result for diversities.

\begin{thm}
For any diversity $(X,\delta)$, the tight span $(T_X,\delta_T)$ is hyperconvex. \label{lem:tightspanConvex}
\end{thm}
\begin{pf}
Let $r \colon \Pf(T_X) \rightarrow \Re$ be given such that for all finite $\sF \subseteq \Pf(T_X)$
\[
\sum_{F \in \sF} r(F) \geq \delta_T \left( \bigcup_{F \in \sF} F \right).
\]
Without loss of generality we can assume $r(\emptyset)=0$.
  We need to find $g \in T_X$ so that $\delta_T(G \cup \{g\})\leq r(G)$ for all $G \subseteq T_X$. 
 
 Define $\omega$ on $\Pf(X)$ by
 \[\omega(A) = \inf_{F \subseteq T_X} \{ r(F) + \Phi_F(A) \colon |F|<\infty  \}.\]
 We have $\omega(\emptyset)=0$.
 Suppose that $\sA \subseteq \Pf(X)$, $|\sA|<\infty$ and let  $\{F_A:A \in \sA\}$ be a collection of finite subsets of $T_X$ indexed by elements of $\sA$. From Lemma~\ref{lem:phi}  and \eqref{eq:deltatPhi} we have
  \begin{eqnarray*}
 \delta\left( \bigcup_{A \in \sA} A \right) & \leq & \delta_T\left( \bigcup_{A \in \sA} F_A \right) + \Phi_{\left(\bigcup_{A \in \sA} F_A\right)} \left(\bigcup_{A \in \sA}A \right) \\
 & \leq & \sum_{A \in \sA} ( r(F_A) + \Phi_{F_A}(A) ),
 \end{eqnarray*} 
 so that
 \[ \delta\left( \bigcup_{A \in \sA} A \right)  \leq  \sum_{A \in \sA} \omega(A),\]
  and $\omega \in P_X$.    
  
  There is $g \in T_X$ such that $g \preceq \omega$. Consider finite $F \subseteq T_X$. Applying Lemma~\ref{lem:deltat}, 
  \begin{eqnarray*}
  \delta_T(F \cup \{g\}) &=& \sup_{A \in \Pf(X)} \left\{g(A) - \Phi_F(A)\right\}  \\
  & \leq & \sup_{A \in \Pf(X)} \left\{ (r(F) + \Phi_F(A)) - \Phi_F(A)\right\}  \\
& = & r(F),
\end{eqnarray*}
as required.
\end{pf}

The metric tight span construction gives an isometric embedding $\kappa$ from a metric space $(X,d)$ into an injective (hyperconvex) metric space. Isbell showed that this embedding is minimal in that no proper subspace of the tight span both contains $\kappa(X)$ and is injective. Such an embedding is called an {\em injective envelope}, and all injective envelopes of a metric space are equivalent \cite[Thm 2.1]{isbell}. 

Here we prove the analogous result for diversities that the embedding $\kappa$ of a diversity into its tight span is also an injective envelope. 

The class of all diversities with all non-expansive maps as morphisms forms a category,
which we will denote $\mathbf{Dvy}$ and call the `Category of Diversities'. The definitions of embeddings and injective objects then correspond to concepts in category theory, as reviewed in \cite{Adamek90}.  Lemma~\ref{lem:kappaisessential}  together with the injectivity of $(T_X,\delta_T)$ establishes that $(T_X,\delta_T)$ is the {\em injective hull} of $(X,\delta)$ in the category $\mathbf{Dvy}$ \cite[pg. 156]{Adamek90}. Proposition 9.20(5) of \cite{Adamek90} demonstrates the equivalence between the category theory injective hull and the  injective envelope introduced in \cite{isbell}.

\begin{lem} \label{lem:kappaisessential}
Let $\phi$ be a non-expanding map from $(T_X,\delta_T)$ to diversity $(Y,\delta_Y)$.  If $\pi = \phi \circ \kappa$ is an embedding from $(X,\delta)$ to $(Y,\delta_Y)$ then $\phi$ is an embedding from $(T_X,\delta_T)$ to $(Y,\delta_Y)$.
\end{lem}
\begin{pf}
Since $\phi$ is non-expanding $\delta_T(F) \geq \delta_Y(\phi(F))$ for all finite $F \subseteq T_X$.  Using part 3 of Theorem~\ref{thm:deltaTmin}  we will show that $\delta_T(F) \leq \delta_Y(\phi(F))$, so that $\phi$ is an embedding. 

Consider $f \in T_X$. 
Define $g$ on $\Pf(X)$ by $g(A) =  \delta_Y(\pi(A) \cup \phi(\{f\}))$ for all finite  $A$. Then
for any finite  $A \subseteq X$ we have
\[g(A) =  \delta_Y(\pi(A) \cup \phi(\{f\}) )= \delta_Y(\phi(\kappa(A) \cup \{f\})) \leq \delta_T(\kappa(A) \cup \{f\}) = f(A)\]
 for all $A$.  For all finite collections $\sA \subseteq \Pf(X)$ we have
\begin{eqnarray*}
\sum_{A \in \sA} g(A) &=& \sum_{A \in \sA}  \delta_Y(\pi(A) \cup \phi(\{f\})) \\
& \geq & \delta_Y\left(\bigcup_{A \in \sA} \pi(A) \right) \\&=& \delta \left(\bigcup_{A \in \sA} A \right),
\end{eqnarray*}
so that $g \in P_X$ and $g \preceq f$. Hence $g(A)=f(A)$ for all finite $A \subseteq X$. It follows that
\[\delta_Y(\pi(A) \cup \phi(\{f\})) =  \delta_T(\kappa(A) \cup \{f\})\]
for all $f \in T_X$ and finite $A \subseteq X$. 

Define $\hdelta$ on $T_X$ by $\hdelta(F) = \delta_Y(\phi(F))$. Then $\hdelta$ is a diversity and $\hdelta(\kappa(Y) \cup \{f\}) = f(Y)$ for all finite $Y \subseteq X$. By Theorem~\ref{thm:deltaTmin}, $\hdelta(F) \geq \delta_T(F)$ for all finite $F$. 

As $\phi$ is non-expansive $\delta_T(F) = \delta_Y(\phi(F))$ for all finite  $F$ and $\phi$ is an embedding.
\end{pf}

The following theorem is a translation of \cite[Proposition 9.20(4)]{Adamek90} to diversities.

\begin{thm}\label{thm:envelope}
If there is an embedding $\pi$ from $(X,\delta)$ into  $(Y,\delta_Y)$  and $(Y,\delta_Y)$ is injective (hyperconvex) then there is an embedding $\phi$ from $(T_X,\delta_T)$ into $(Y,\delta_Y)$ such that $\pi = \phi \circ \kappa$.
\end{thm}
\begin{pf}
Since $\pi$ is a non-expansive map, $(Y,\delta_Y)$ is injective, and $\kappa$ is an embedding of $(X,\delta)$ into $(T_X,\delta_T)$, there is a non-expansive map $\phi:T_X \rightarrow Y$ such that $\pi = \phi \circ \kappa$.   By Lemma~\ref{lem:kappaisessential}, $\phi$ is an embedding.
\end{pf}

\begin{cor}
Let $(X,\delta)$ be a diversity. The following are equivalent:
\begin{enumerate}
\item $(X,\delta)$ is hyperconvex;
\item $(X,\delta)$ is injective;
\item There is an isomorphism between $(X,\delta)$ and its tight span, $(T_X,\delta_T)$.
\end{enumerate}
\end{cor}

\begin{pf}
Parts 1 and 2 are equivalent by Theorem \ref{lem:injecthyper}.  To see that part 2 implies part 3, let $(Y,\delta_Y) = (X,\delta)$ and $\pi = \mbox{id}$ in Theorem~\ref{thm:envelope}.    Then there is an embedding $\phi$ from $(T_X,\delta_T)$ to $(X,\delta)$ such that $\phi \circ \kappa = \mbox{id}$.  So $\kappa$ is surjective and part 3 follows.
 Finally,  since hyperconvexity is invariant under isomorphism, part 1 follows from part 3.
\end{pf}

\section{Tight span of the diameter diversity} \label{sec:tightlinks}

In this section we prove that tight-span theory for metrics is embedded within the tight-span theory for diversities. The link between the two is provided by the {\em diameter diversity} as introduced above.

\begin{defn}
Given a metric space $(X,d)$ we define the function $\delta = \diam_d$ by 
\[\delta(A) = \diam_d(A) = \max \{d(a,a'):a,a' \in A \}\]
 for finite $A \subseteq X$, with $\diam_d(\emptyset) = 0$. We call $(X,\diam_d)$ the {\em diameter diversity} for $(X,d)$. 
\end{defn}

Note that if we restrict $\diam_d$ to pairs of elements we recover $d$ as the induced metric.  We will establish close links between tight spans of metrics and tight spans of their diameter diversities. 

\[\begin{CD}
(X,d)  	@>\mbox{\tiny tight span}>>	(T^d_X,d_T) \\
@V{\delta=\diam_d}VV			@VV{\delta_T = \diam_{d_{\footnotesize T}}}V \\
(X,\delta) @>\mbox{\tiny tight span}>> (T^\delta_X,\delta_T)\\
\end{CD}\]
\vspace{0.5cm}

\begin{lem} \label{lem:diamInject} 
\begin{enumerate} 
\item Let $(Y,\delta)$ be a diversity with induced metric $(Y,d_\delta)$. Let $(X,d)$ be a metric space and let $(X,\diam_d)$ be the associated diameter diversity. Then $\phi$ is a non-expansive map from $(Y,\delta)$ to $(X,\diam_d)$ if and only if it is a non-expansive map from $(Y,d_\delta)$ to $(X,d)$.
\item A metric space $(X,d)$ is injective (hyperconvex) if and only if the diameter diversity $(X,\diam_d)$ is injective (hyperconvex).
\item The tight span $(T^\delta_X,\delta_T)$ of a diameter diversity is itself a diameter diversity. \end{enumerate}
\end{lem}
\begin{pf}
1. Suppose that $\phi$ is a non-expansive map from $(Y,\delta)$ to $(X,\diam_d)$. For all $y_1,y_2 \in Y$ we have
\[d_\delta(y_1,y_2) = \delta(\{y_1,y_2\}) \geq \diam_d(\{\phi(y_1),\phi(y_2)\}) = d(\phi(y_1),\phi(y_2)) ,\]
so $\phi$ is non-expansive from $(Y,d_\delta)$ to $(X,d)$. Conversely, suppose  $\phi$ is a non-expansive map from  $(Y,d_\delta)$ to $(X,d)$. Then for any finite $A \subseteq Y$ we have
\begin{eqnarray*}
\delta(A)  &\geq &  \sup\{d_\delta(a_1,a_2):a_1,a_2 \in A \} \\&  \geq& \sup\{d(\phi(a_1),\phi(a_2)):a_1,a_2 \in A\}\\& =& \diam_d(\phi(A)).\end{eqnarray*}
2. Suppose that $(X,d)$ is injective. Let $(Y_1,\delta_1)$, $(Y_2,\delta_2)$ be two diversities with induced metrics $d_1,d_2$. Let $\pi$ be an embedding from $(Y_1,\delta_1)$ into $(Y_2,\delta_2)$ and let $\phi$ be a non-expansive map from $(Y_1,\delta_1)$ to $(X,\diam_d)$. Then $\pi$ embeds $(Y_1,d_1)$ into $(Y_2,d_2)$, and by part 1, $\phi$ is a non-expansive map from $(Y_1,d_1)$ to $(X,d)$. As $(X,d)$ is an injective metric space there is a non-expansive map $\psi$ from $(Y_2,d_2)$ to $(X,d)$ such that $\phi = \psi \circ \pi$, which by part 1 is a non-expansive map from $(Y_2,\diam_{d_2})$ to $(X,\diam_d)$.  Since $\delta_2(A) \geq \diam_{d_2}(A)$ for all $A$, $\psi$ is non-expansive from $(Y_2,\delta_2)$ to $(X,\diam_d)$. Hence $(X,\diam_d)$ is injective.

Conversely, suppose $(X,\diam_d)$ is an injective diversity. Let  $(Y_1,d_1)$, $(Y_2,d_2)$ be two metric spaces, let $\pi$ be an embedding of $(Y_1,d_1)$ into $(Y_2,d_2)$, and let $\phi$ be a non-expansive map from $(Y_1,d_1)$ to $(X,d)$. Then $\phi$ is a non-expansive map from $(Y_1,\diam_{d_1})$ to $(X,\diam_d)$ and since $(X,\diam_d)$ is injective, there is a non-expansive map $\psi$ from $(Y_2,\diam_{d_2})$ to $(X,\diam_d)$ such that $\phi = \psi \circ \pi$. Applying part 1 again, we have that $\psi$ is the required non-expansive map from $(Y_2,d_2)$ to $(X,d)$. Hence $(X,d)$ is injective. \\
3. Since $(X,\delta)$ is a diameter diversity, for any finite $F$ and $\{ A_f \}_{f \in F} \subseteq \Pf(X)$, we have
\[
\delta\left(\bigcup_{f \in F} A_f \right)= \delta(A_{f_1} \cup A_{f_2})
\]
for some $f_1, f_2 \in F$.
 Hence for finite $F \subseteq T_X^\delta$
\begin{eqnarray*}
\delta_T(F) & = & \sup_{A_f} \left\{ \deltap{\bigcup_{f \in F} A_f} - \sum_{f \in F} f(A_f) \right\} \\
& = & \max_{f_1,f_2 \in F} \sup_{A_1,A_2 \in \Pf(X)} \left\{ \delta(A_1 \cup A_2) - f_1(A_1) - f_2(A_2) \right\} \\
& = & \max_{f_1,f_2 \in F} \delta_T(\{f_1,f_2\}).
\end{eqnarray*}
\end{pf}

\begin{thm} \label{thm:metricdiamdiverseequiv}
Let $(X,d)$ be a metric space with metric tight span $(T^d_X,d_T)$.
Let $(X, \delta)$ be the associated diameter diversity where $\delta = \diam_d$, and let $(T^\delta_X,\delta_T)$ be its diversity tight span.  Then
\begin{enumerate}
\item The metric space obtained by restricting $\delta_T$ to pairs in $T^\delta_X$ is  isometric to the metric space $(T^d_X,d_T)$.   
\item  The diversity obtained by taking the diameter on the metric space $(T^d_X,d_T)$ is isomorphic to the diversity $(T^\delta_X,\delta_T)$.
\end{enumerate}
\end{thm}
\begin{pf}
First note that for any metric spaces $(X_1,d_1)$ and $(X_2,d_2)$  a map $\phi$ from $X_1$ to $X_2$ is an embedding from $(X_1,d_1)$ to $(X_2,d_2)$ if and only if $\phi$ is an embedding from $(X_1,\diam_{d_1})$ to $(X_2,\diam_{d_2})$.

Let $(T^d_X,\delta_{d_T})$ be the diameter diversity associated to $(T^d_X,d_T)$ and let $(T^\delta_X,d_{\delta_T})$ be the induced metric for $(T^\delta_X,\delta_T)$.  Let $\kappa_d$ be the Kuratowski embedding from $(X,d)$ to $(T^d_X,d_T)$. Then $\kappa_d$ is also an embedding from $(X,\delta)$ to $(T^d_X,\delta_{d_T})$. In the same way, let $\kappa_\delta$ be the Kuratowski embedding from $(X,\delta)$ to $(T^\delta_X,\delta_T)$. Then $\kappa_\delta$ is also an embedding from $(X,d)$ to $(T^\delta_X,d_{\delta_T})$.

By Lemma~\ref{lem:diamInject} 2.,  $(T^d_X,\delta_{d_T})$ is a hyperconvex diversity and $(T^\delta_X,d_{\delta_T})$ is a hyperconvex metric space. Applying 
\cite[Proposition 9.20(4))]{Adamek90}  in the category {\bf Met}  there is an embedding $\phi$ from $(T^d_X,d_T)$ to $(T^\delta_X,d_{\delta_T})$ such that 
\begin{equation} \kappa_\delta = \phi \circ \kappa_d. \label{eq:diam1} \end{equation} 

The identity map $\mathrm{id}_{T^d_X}$ on $(T^d_X,d_T)$ is non-expansive and $\phi$ is an embedding, so applying the definition of injective metric spaces to $(T^d_X,d_T)$ we have that there is a non-expansive map $\psi$ from
$T^\delta_X$ to $T^d_X$ such that 
\begin{equation}
\psi \circ \phi = \mathrm{id}_{T^d_X}. \label{eq:diam2}
\end{equation}

By Lemma~\ref{lem:diamInject} 3., the diversity $(T_X^\delta,\delta_T)$ is a diameter diversity and so from Lemma~\ref{lem:diamInject} {1.}, the map $\psi$ is also a non-expansive map from $(T_X^\delta,\delta_T)$ to $(T^d_X,\delta_{d_T})$. Combining \eqref{eq:diam1} and \eqref{eq:diam2} we have
\begin{align*}
\psi \circ \kappa_\delta & = \psi \circ \phi \circ \kappa_d \\
& = \mathrm{id}_{T^d_X} \circ \kappa_d 
\end{align*}
which is an embedding. By Lemma~\ref{lem:kappaisessential} we have that  $\psi$ is an embedding, implying that that
$\phi$ is both an  isomorphism from $(T^d_X,d_T)$ to $(T^\delta_X,d_{\delta_T})$ and an isomorphism from $(T^d_X,\delta_{d_T})$ to $(T^\delta_X,\delta_T)$.
\end{pf}

\section{Phylogenetic diversity} \label{sec:phylo}

\newcommand{\tree}{t}
\newcommand{\dtree}{\delta_{\tree}}

A metric space $(X,d)$ is {\em additive} or {\em tree-like} if there is a tree with nodes partially labelled by $X$ so that for each $x,y \in X$ the length of the path (including branch-lengths) connecting $x$ and $y$ equals $d(x,y)$. Dress \cite{Dress84} showed that if $(X,d)$ is additive then its metric tight span corresponds exactly to the smallest tree it can be embedded in. The elements of the tight span correspond not only to the nodes of the original tree, but also the points along the edges. Here we will prove analogous results about phylogenetic diversity.

Following \cite{Dress84} we will work with {\em $\Re$-trees}  (also called metric-trees), rather than graph-theoretic trees.
\begin{defn}\cite{espinola_kirk,chiswell}
\begin{enumerate}
\item Let $(\sX,d)$ be a metric space and let $x,y$ be two points at distance $d(x,y) = r$. A \emph{geodesic} joining $x,y$ is a map $c:[0,r] \rightarrow \sX$ such that $c(0) = x$, $c(r) = y$ and $d(c(s),c(t)) = |t-s|$ for all $s,t \in [0,r]$. The image of $c$ is called a {\em geodesic segment}.
\item \cite[Defn 2.1]{espinola_kirk}
A  metric space $(\sX,d)$ is an $\Re$-tree if
\begin{enumerate}
\item there is a unique geodesic segment $[x,y]$ joining each pair of points $x,y \in \sX$.
\item if $[y,x] \cap [x,z] = \{x\}$ then $[y,x] \cup [x,z] = [y,z]$.
\end{enumerate} 
Hence if $x,y,z$ are three points in an $\Re$-tree then 
\begin{equation}
[x,y] \subseteq [x,z] \cup [y,z]. \label{eq:rt}
\end{equation}
\end{enumerate}
\end{defn}

Phylogenetic diversity, as introduced by \cite{Faith92} and investigated extensively by \cite{Steel05,Minh09,Pachter04} and others, can be viewed as a generalisation of additive metrics. The phylogenetic diversity of a set of nodes or points in a tree is the length of the smallest subtree connecting them, so that the restriction of a phylogenetic diversity to pairs of points gives an additive metric. A formal definition of phylogenetic diversity on $\Re$-trees requires a bit more machinery. 

For a  $\Re$-tree $(\sX,d)$, let $\mu$ be the one-dimensional Hausdorff  measure on it \cite{Edgar08}.   The important features of $\mu$ for our purposes is that it is defined on all Borel sets, it is monotone, and it is additive on disjoint sets. 
 Furthermore, for any points $a,b \in \sX$, $\mu([a,b])=d(a,b)$, and naturally $\mu(\{a\})=0$. 
See \cite{evans_winter} for a related measure on $\Re$-trees.

\begin{defn}
\begin{enumerate}
\item The convex hull of a set $A \subseteq \sX$ is
\[\conv{A} = \bigcup_{a,b \in A} [a,b]\]
and we say that $A$ is convex if $A = \conv{A}$.
\item Let $(\sX,d)$ be an $\Re$-tree. The {\em real-tree diversity} $(\sX,\delta_\tree)$ for $(\sX,d)$ is defined by
\[
\delta_\tree(A) := \mu(\conv{A})
\]
for all finite $A \subseteq \sX$. Note that since $A$ is finite, $\conv{A}$ is closed and hence $\mu(\conv{A})$ is defined.
\end{enumerate}
\end{defn}

First we prove that this phylogenetic diversity satisfies the diversity axioms (D1) and (D2).

\begin{thm} \label{thm:phyloIsDiv}
Let $(\sX,d)$ be an $\Re$-tree. Then $(\sX,\delta_\tree)$ is a diversity.
\end{thm}
\begin{pf}
Since $\mu$ is a measure, $\delta_\tree$ is non-negative and also monotonic. If $|A| \leq 1$ then $\conv{A}=A$  and so $\delta_\tree(A) = \mu(A)=0$. If $|A|>1$ then select distinct $a,b\in A$.  Since $\conv{[a,b]}=[a,b]$ and $\mu([a,b])= d(a,b)$ we have $\delta_\tree(A) \geq \delta_\tree(\{a,b\})= d(a,b)>0$. This proves (D1).

Let $A,B,C \in \Pf(\sX)$ and suppose that $B \neq \emptyset$. From \eqref{eq:rt} we have
\begin{equation} \label{eqn:geodesiclemma}
[a,c] \subseteq [a,b] \cup [b,c]
\end{equation}
for all $a \in A$, $b \in B$ and $c \in C$. Hence
\[\conv{A \cup C} \subseteq \conv{A \cup B} \cup \conv{B \cup C}\]
and
\begin{eqnarray*}
\delta_\tree(A \cup C) & = & \mu(\conv{A \cup C}) \\
&\leq& \mu(\conv{A \cup B}) + \mu(\conv{B \cup C}) \\
& = & \delta_\tree(A \cup B) + \delta_\tree(B \cup C),
\end{eqnarray*}
giving us the triangle equality (D2).  

\end{pf}

%
%
%
%
%
%
%

We now show that complete real-tree diversities are hyperconvex, proving the diversity analogue of \cite[Theorem 3.2]{Kirk98}.

\begin{lem} \label{lem:closedconvex}
Let $(\sX,d)$ be an $\Re$-tree with associated tree diversity $(\sX,\delta_\tree)$. For all finite $C \subseteq \sX$ and $r \geq \delta_\tree(C)$, the ball $B(C,r) = \{x \in \sX:\delta_t(C \cup \{x\}) \leq r\}$ is closed and convex.
\end{lem}
\begin{pf}
For any finite but non-empty $C \subseteq \sX$ the function 
\[\phi: \sX \rightarrow \Re : x \mapsto \delta_\tree(C \cup \{x\} )\]
 is continuous. Hence when $r \geq \dtree(C)$ the ball
\[ B(C,r) := \phi^{-1}(A) =  \{x \in T_X^d :\delta_\tree(C \cup \{x\}) \leq r\}\]
is closed. 

To prove convexity, suppose that  $x_1,x_2 \in B(C,r)$. Fix $a \in C$.  For all $y \in [a,x_1]$, $\conv{C \cup \{y\} } \subseteq \conv{ C \cup \{x_1\} }$ and so $\delta_t( C \cup \{y\}  ) \leq \delta_t(  C \cup \{x_1\} )$ showing that $y \in B(C,r)$.  We have that $[a,x_1]$, and by symmetry $[a,x_2]$, are contained in $B(C,r)$.
By \eqref{eq:rt} we have 
\[[x_1,x_2] \subseteq [a,x_1] \cup [a,x_2] \subseteq B(C,r)\]
so that $B(C,r)$ is both closed and convex.
\end{pf}

\begin{thm} \label{thm:phyloIsHyper}
Let $(\sX,d)$ be an $\Re$-tree with associated real-tree diversity $(\sX,\delta_\tree)$. Then $(\sX,\delta_\tree)$ is hyperconvex if and only if $(\sX,d)$ is complete. 
\end{thm}
\begin{pf}

Suppose that $(\sX,d)$ is a complete $\Re$-tree. Then $(\sX,d)$ is a hyperconvex metric space \cite[Theorem 3.2]{Kirk98}. Suppose that $r:\Pf(\sX) \rightarrow \Re$ satisfies 
\[\delta_\tree \left(\bigcup_{A \in \sA} A \right) \leq \sum_{A \in \sA} r(A)\]
for all finite $\sA \subseteq \Pf(\sX)$. We will show that the collection of balls 
\[\Gamma = \{B(A,r(A)):A \in \Pf(\sX)\}\]
has a non-empty intersection. 

First we show that the members of $\Gamma$ intersect pairwise.
Consider a pair of nonempty finite subsets  $A_i,A_j$ of $X$.  To show that $B(A_i,r(A_i))$ and $B(A_j,r(A_j))$ intersect, we show that there is $v$ such that $\dtree(A_i \cup \{v\}) \leq r(A_i)$ and $\dtree(A_j \cup \{v\}) \leq r(A_j)$. This clearly holds if there is $\conv{A_i} \cap \conv{A_j} \neq \emptyset$. Suppose then that $\conv{A_i}$ and $\conv{A_j}$ are disjoint.  Since $A_i, A_j$ are finite, $\conv{A_i}$ and $\conv{A_j}$ are closed subtrees of $T_X^d$. By \cite[Ch. 2, Lemma 1.9]{chiswell} there exists $a_i \in \conv{A_i}$ and $a_j \in \conv{A_j}$ such that $[a_i,a_j] \cap \conv{A_i} = \{a_i\}$ and $[a_i,a_j] \cap \conv{A_j} = \{a_j\}$ and for all $x \in A_i$ and $y \in A_j$ we have $[a_i,a_j] \subseteq [x,y]$. 
Then, 
\begin{eqnarray*}
r(A_i) + r(A_j) & \geq & \delta_t(A_i \cup A_j) \\
  & = & \mu( \conv{A_i \cup A_j) } \\ 
 & \geq & \mu(\conv{A_i}) + \mu([a_i,a_j]) + \mu(\conv{A_j}) \\
& = & \dtree(A_i) + d(a_i,a_j) + \dtree(A_j).
\end{eqnarray*}
Hence there is $v \in [a_i,a_j]$ such that $d(a_i,v) \leq r(A_i) - \dtree(A_i)$ and $d(a_j,v) \leq r(A_j) - \dtree(A_j)$, so that
\begin{eqnarray*}
\dtree(A_i \cup \{v\}) & = & \dtree(A_i ) + \dtree(\{a_i,v\}) \\
& = & \delta(A_i) + d(a_i,v) \\
& \leq & r(A_i),
\end{eqnarray*}
and likewise $\dtree(A_j \cup \{v\}) \leq r(A_j)$.  

We have established that $\Gamma$ satisfies the pairwise intersection property. The closed, convex sets of an $\Re$-tree satisfy the Helly property \cite{espinola_kirk},  so every finite subcollection of $\Gamma$ has non-empty intersection. By the completeness of $(\sX,d)$,  $\Gamma$ has a non-empty intersection, so there is $v$ such that  $\delta_\tree(A \cup \{v\}) \leq r(A)$ for all $A \in \Pf(\sX)$. This proves that $(\sX,\delta_\tree)$ is hyperconvex.

For the converse, we note that completeness of $(\sX,d)$ follows directly from \cite[Proposition 3,2]{Espinola01} and the 
definition of hyperconvexity for diversities.
\end{pf}

\begin{defn}
A diversity $(X,\delta)$ is a {\em phylogenetic diversity}  if it can be embedded in  a real-tree diversity $(\sX,\delta_\tree)$ for some complete $\Re$-tree $(\sX,d)$. 
\end{defn}

Clearly, every real-tree diversity is a phylogenetic diversity, but a phylogenetic diversity is a real-tree diversity only if its induced metric is an $\Re$-tree.

\begin{thm} \label{thm:phylo_if_tight_is_tree}
Let $(X,\delta)$ be a diversity. Then  $(X,\delta)$ is a phylogenetic diversity if and only if  $(T_X,\delta_T)$ is a real-tree diversity.
\end{thm}
\begin{pf}
Since $(X,\delta)$ is a phylogenetic diversity there is a complete $\Re$-tree $(\sX,d)$ with real-tree diversity $(\sX,\delta_\tree)$ for which there is an embedding $\phi$ from $(\sX,\delta)$ into $(\sX,\delta_\tree)$. By Theorem~\ref{thm:phyloIsHyper} $(\sX,\delta_\tree)$ is hyperconvex. By Theorem~\ref{thm:envelope} there is an embedding $\psi$ from $(T_X,\delta_T)$ into $(\sX,\delta_\tree)$ such that $\phi = \psi \circ \kappa$.

Let $(T_X,d_{\delta_T})$ be the induced metric for $(T_X,\delta_T)$. It follows directly from the hyperconvexity of $(T_X,\delta_T)$ that $(T_X,d_{\delta_T})$ is convex. For any $f,g \in T_X$ and geodesic segment $[f,g]$ in $T_X$, the image of $[f,g]$ under $\psi$ is the unique geodesic segment between $\psi(f)$ and $\psi(g)$. It follows that $\psi(T_X)$ is a convex subset of $(\sX,d)$ and $(\sX,d)$ restricted to $\psi(T_X)$ is an $\Re$-tree \cite[pg. 36]{chiswell}. Restricting $(X,\delta_t)$ to $\psi(T_X)$ then gives a real-tree diversity which is isomorphic to $(T_X,\delta_T)$.

For the converse, note that the map $\kappa$ from $(X,\delta)$ into its tight span is an embedding, so that $(X,\delta)$ is a phylogenetic diversity.
\end{pf}

We now link the $\Re$-tree given by the diversity tight span of a phylogenetic diversity and the tight span of its induced metric. 

\begin{lem} \label{lem:divembed}
Let $(X,d_X)$ and $(Y,d_Y)$ be complete $\Re$-trees and let $(X,\delta_X)$ and $(Y,\delta_Y)$ be the associated real-tree diversities. Then 
\begin{enumerate}
\item$\psi:X \rightarrow Y$ is a non-expansive map from $(X,d_X)$ to $(Y,d_Y)$ if and only if it is a non-expansive map from $(X,\delta_X)$ to $(Y,\delta_Y)$.
\item $\psi:X \rightarrow Y$ is an embedding from $(X,d_X)$ to $(Y,d_Y)$ if and only if it is an embedding from $(X,\delta_X)$ to $(Y,\delta_Y)$.
\end{enumerate}
\end{lem}
\begin{pf}
1. Suppose that $\psi \colon X \rightarrow Y$ is a non-expansive map from $(X,d_X)$ to $(Y,d_Y)$.  For any finite $A \subset X$, we have $\delta_Y( \phi(A)) = \mu( \conv{ \phi(A) }$.
First note that $\conv{\phi(A)} = \phi( \conv{A})$.  Then note that since in this case the one-dimensional Hausdorff measure of a set is a limit of infima of the total length of countable covers of a set by geodesic segments \cite[Section 6.1]{Edgar08}, $\mu(\phi(B)) \leq \mu(B)$ for all measurable $B \subseteq X$.  This proves that $\psi$ is a non-expanding map with respect to diversities.  The other direction is immediate. 

\noindent 2.  The argument follows as in part 1, with showing $\mu(\phi(B)) = \mu(B)$ for all measurable $B \subseteq X$. 
\end{pf}

%
%

\begin{thm}
Let $(X,\delta)$ be a phylogenetic diversity and let $(X,d)$ be its induced metric. Let $(T^\delta_X,\delta_T)$ be the diversity tight span of $(X,\delta)$ and let $(T^d_X,d_T)$ be the metric tight span of $(X,d)$. Then $(T^d_X,d_T)$  is isometric with the induced metric of $(T^\delta_X,\delta_T)$.
\end{thm}
\begin{pf}

By \cite[Theorem 8]{Dress84}, $(T^d_X,d_T)$ is an $\Re$-tree. Let $(T^d_X,\delta_{d_T})$ be the corresponding real-tree diversity, which is hyperconvex by Theorem~\ref{thm:phyloIsHyper}. Let $(T_X^\delta,d_\delta)$ denote the induced metric of $(T^\delta_X,\delta_T)$. From Theorem~\ref{thm:phyloIsHyper} we have that $(T_X^\delta,d_\delta)$ is a complete $\Re$-tree and is therefore a hyperconvex metric space \cite[Theorem 3.2]{Kirk98}.

Let $\kappa_d$ be the Kuratowski embedding from $(X,d)$ to $(T^d_X,d_T)$ and let $\kappa_\delta$ be the Kuratowski embedding from $(X,\delta)$ to $(T^\delta_X,\delta_T)$. The map $\kappa_\delta$ is then also an embedding between the induced metric $(X,d)$ and the induced metric $(T^\delta_X,d_{\delta_T})$. Applying 
\cite[Proposition 9.20(4))]{Adamek90}  in the category {\bf Met}, there is an embedding $\phi:(T^d_X,d_T) \rightarrow (T^\delta_X,d_{\delta_T})$ such that 
\begin{equation} \kappa_\delta = \phi \circ \kappa_d. \label{eq:div1} \end{equation} 
By Lemma~\ref{lem:divembed}, $\phi$ is also an embedding from the diversity $(T_X^d,\delta_{d_T})$ to the diversity $(T_X^\delta,\delta_T)$. For all $A \in \Pf(X)$,
\[\delta(A) = \delta_T(\kappa_\delta(A)) = \delta_T(\phi(\kappa_d(A))) = \delta_{d_T}(\kappa_d(A))\]
so that $\kappa_d$ embeds $(X,\delta)$ in $(T_X^d,\delta_{d_T})$.

The identity map $\mathrm{id}_{T^d_X}$ on $(T^d_X,d_T)$ is non-expansive and $\phi$ is an embedding, so applying the definition of injective metric spaces to $(T^d_X,d_T)$ we have that there is a non-expansive map $\psi$ from
$(T^\delta_X,d_{\delta_T})$ to $(T^d_X,d_T)$ such that 
\begin{equation}
\psi \circ \phi = \mathrm{id}_{T^d_X}. \label{eq:div2}
\end{equation}

Applying part 1 of  Lemma~\ref{lem:divembed} we see that the map $\psi$ is also a non-expansive map from $(T_X^\delta,\delta_T)$ to $(T^d_X,\delta_{d_T})$. Combining \eqref{eq:div1} and \eqref{eq:div2} we have
\begin{align*}
\psi \circ \kappa_\delta & = \psi \circ \phi \circ \kappa_d \\
& = \mathrm{id}_{T^d_X} \circ \kappa_d 
\end{align*}
which is an embedding. By Lemma~\ref{lem:kappaisessential} we have that  $\psi$ is an embedding, implying that that
$\phi$ is  an  isometry from $(T^d_X,d_T)$ to $(T^\delta_X,d_{\delta_T})$.
\end{pf}

\section{Tight span and the Steiner tree problem} \label{sec:Steiner}

Let $X$ be a finite set of points
in a metric space $(M,d)$. The {\em (metric) Steiner tree problem} is to find the shortest network that connects them. Clearly this network will always be a tree. More formally

\vspace{0.1in}

\noindent {\sc Metric Steiner  Problem.} \\
\noindent {\bf Input}: Subset $X$ of a metric space $(M,d)$.\\
\noindent {\bf Problem:} Find a (graph theoretic) tree $T$  for which $X \subseteq V(T) \subseteq M$ and
\[\sum_{\{u,v\}\in E(T)} d(u,v)\]
is minimised.\\

Dress and Kr\"uger \cite{Dress87} examined an `abstract' metric Steiner problem where one  drops the constraint that $V(T) \subseteq M$. This abstract Steiner tree was  one of the first distance-based criteria proposed for the inference of phylogenetic trees \cite{Beyer74,Waterman77}, though it is now not widely used.  Suppose that $T$ is a tree with edge weights $w:E(T)\rightarrow \Re_{\geq 0}$. Given $u,v \in V(T)$ we let $d_w(u,v)$ denote the sum of edge weights along the path from $u$ to $v$. \\

\noindent {\sc Abstract Steiner  Problem.} \\
\noindent {\bf Input}: Finite metric space  $(X,d)$.\\
\noindent {\bf Problem:} Find a (graph theoretic) tree $T$ and edge weighting $w:E(T) \rightarrow \Re$ such that $X \subseteq V(T)$, $d_w(x,y) \geq d(x,y)$ for all $x,y \in X$ and 
\[\sum_{e\in E(T)} w(e)\]
is minimised.\\

Suppose that $T$ is a solution to the metric Steiner problem for $X \subseteq M$. Define the weight function $w:E(T) \rightarrow \Re$ by $w(\{u,v\}) = d(u,v)$. Then, by the triangle inequality, $d_w(x,y) \geq d(x,y)$ for all $x,y \in X$. It follows then that the length of the minimum abstract Steiner tree for $(X,d|_X)$ is a lower bound for the metric Steiner problem. Dress and Kr\"uger showed that the lower bound becomes tight when $(M,d)$ equals $(T_X,d_T)$, the metric tight span of $X$.

\begin{thm}[\cite{Dress87}] \label{thm:DKTheorem}
Let $(X,d)$ be a finite metric space. For every solution $(T,w)$ to the abstract Steiner tree problem there is a map $\phi:V(T) \rightarrow T_X$ such that $\phi(x) = \kappa(x)$ for all $x \in X$ and $w(\{u,v\}) = d_T(\phi(u),\phi(v))$ for all $\{u,v\} \in E(T)$.
\end{thm}
 Hence the length of the minimal Steiner tree for $\kappa(X)$ in $(T_X,d_T)$ equals the length of the minimal abstract Steiner tree for $(X,d)$ and the minimal abstract Steiner trees can be embedded within the tight span. A direct corollary is that if $d$ is tree-like then the abstract Steiner tree equals the tree corresponding to $d$.

Here we show that, using diversities, we can obtain a tighter bound on the metric Steiner problem than that given by the abstract Steiner problem.
 Given a tree $T$ with edge weights $w$ and $A \subseteq V(T)$ we let $\delta_w(A)$ be the sum of edge weights in the smallest subtree of $T$ connecting $A$. Hence $(X,\delta_w|_X)$ is a phylogenetic diversity.\\

\noindent {\sc Diversity Steiner  Problem.} \\
\noindent {\bf Input}: Finite diversity  $(X,\delta)$.\\
\noindent {\bf Problem:} Find a (graph theoretic) tree $T$ and edge weighting $w:E(T) \rightarrow \Re$ such that $X \subseteq V(T)$, $\delta_w(Y) \geq \delta(Y)$ for all $Y \subseteq X$, and
\[\sum_{e\in E(T)} w(e)\]
is minimised.\\

Let $X$ be a finite subset of a metric space $(M,d)$. For each  $A \subseteq X$ let $\ell(A)$ denote the minimum length of a (metric) Steiner tree connecting the points $A$ in the metric space $(M,d)$. 
We see that $(X,\ell)$ is a diversity.
For each $k\geq 2$, consider the {\em truncated} diversity $\delta^{(k)}$ defined by
\[\delta^{(k)}(A) = \max \{ \ell(B):|B| \leq k,\,\,B \subseteq A \}\]
for all $A \subseteq X$. 

\begin{prop} \label{prop:treebound}
If $(T,w)$ is a minimum length solution for the diversity Steiner problem applied to $\delta^{(k)}$ then the length $\sum_{e \in E(T)} w(e)$ of $T$ is a lower bound for $\ell(X)$, the optimal length of a metric Steiner tree for $X$.
\end{prop}
\begin{pf}
Let $(T',w')$ be a solution to the metric Steiner problem and let $\delta_{w'}$ be the associated phylogenetic diversity. Then for all $B$ such that $|B| \leq k$ we have that $\delta_{w'}(B)$, the length of $T'$ restricted to $B$, is bounded below by $\ell(B) = \delta^{(k)}(B)$. It follows that $\delta^{(k)}(A) \leq \delta_{w'}(A)$ for all $A \subseteq X$, so that $(T',w')$ is a potential solution for the diversity Steiner problem. As $(T,w)$ is optimal, we have
\[\sum_{e \in E(T)} w(e) \leq \sum_{e \in E(T')} w'(e) = \ell(X).\]
\end{pf}

For $k=2$, the bounds provided by the Proposition \ref{prop:treebound}  coincide with those given by length of the minimum abstract Steiner tree.  As $k$ increases, the bounds 
returned by the diversity Steiner tree applied to $\delta^{(k)}$ will tighten, until eventually the diversity Steiner tree will coincide with the metric Steiner tree. Furthermore, we have a direct extension of Theorem~\ref{thm:DKTheorem}, stating that these diversity Steiner trees will all be contained in the diversity tight span.

\begin{thm}
Let $(X,\delta)$ be a finite diversity. For every solution $(T,w)$ to the diversity Steiner tree problem for $(X,\delta)$ there is a map $\phi:V(T) \rightarrow T_X$ such that $\phi(x) = \kappa(x)$ for all $x \in X$ and $w(\{u,v\}) = \delta_T(\{\phi(u),\phi(v)\})$ for all $\{u,v\} \in E(T)$.
\end{thm}
\begin{pf}
Let $\delta_w$ be the diversity on $V(T)$ given by $(T,w)$,
as defined above.
Since $(T,w)$ solves the diversity Steiner problem, $\delta_w(A) \geq \delta(A)$ for all $A \subseteq X$. Let $\kappa$ denote the canonical embedding from $X$ to $T_X$. Then $\kappa$ is a non-expansive map from $(X,\delta_w|_X)$ to $(T_X,\delta_T)$.

The tight span $(T_X,\delta_T)$ is injective. Hence there is a non-expansive map $\phi$ from $(V(T),\delta_w)$ to $(T_X,\delta_T)$ such that $\phi(x) = \kappa(x)$ for all $x \in X$. For each $u,v$ let $w'(\{u,v\}) = \delta_T(\{\phi(u),\phi(v)\})$. Then
\[
w(\{u,v\}) = \delta_w(\{u,v\}) \geq \delta_T(\{\phi(u),\phi(v)\}) = w'(\{u,v\})
\]
for all $u,v \in V$. 

Consider $A \subseteq X$, and let $E_A$ be the set of edges in the smallest subtree of $T$ containing $A$. By the triangle inequality,
\[
\delta_{w'}(A) = \sum_{e \in E_A} w'(e) \geq \delta(X).
\]
Hence $(T,w')$ is a candidate for the diversity Steiner problem, but since $(T,w)$ is already minimum, $\sum_{e \in E(T)} w(e) \leq \sum_{e \in E(T)} w'(e)$. It follows that $w(e)=w'(e)$ for all $e \in E(T)$. 
\end{pf}

\bibliographystyle{elsarticle-num}
\bibliography{tightSpan.bib}

\pagebreak

\begin{figure}[htbp] 
   \centering
   \includegraphics[width=4in]{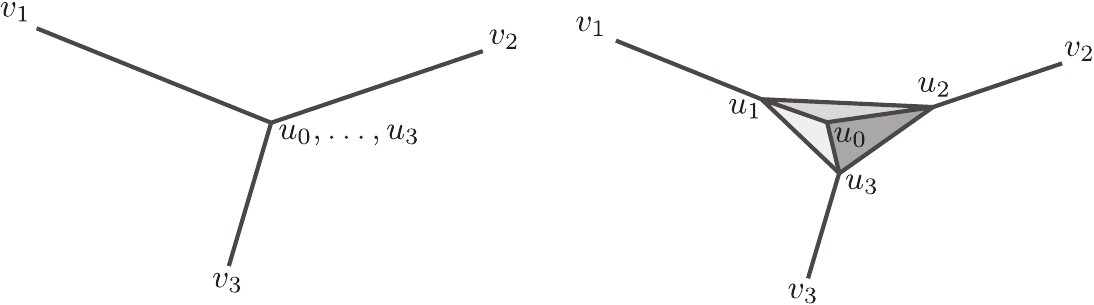} 
  \end{figure}
   
   \vspace*{4cm}
Figure 1.  Two examples of the tight span on three points, with different values for $d(\{1,2,3\})$. On the left an example where $2d_{123} \leq  d_{12} +d_{23}+d_{13}$, and the diversity tight span is one-dimensional and resembles the tight span of the induced metric. On the right a case with $2d_{123} >  d_{12} +d_{23}+d_{13}$, where the diversity consists of a three-cell with three adjacent one-cells.
\pagebreak

\section*{Vitae}
\vspace*{2cm}

\noindent {\bf David Bryant} is Associate Professor of Mathematical Biology at the University of Otago, Dunedin, NZ. His research area is mathematical, computational, and statistical aspects of evolutionary biology. Bryant obtained his Ph.D. in mathematics at the University of Canterbury in 1997 and has since held tenured positions at McGill University and the University of Auckland. \\
\vspace*{1cm}

\noindent  

\noindent {\bf Paul Tupper} is Associate Professor of Mathematics at Simon Fraser University. His research area is the analysis and simulation of models arising in materials science and biology, especially stochastic differential equations, phase field models, and Hamiltonian systems. Tupper obtained his Ph.D. at Stanford University in 2002 and has held positions at McGill University and Simon Fraser University.

\end{document}